\documentclass[conference,letterpaper, onecolumn]{IEEEtran}

\usepackage[utf8]{inputenc} 
\usepackage[T1]{fontenc}
\usepackage{url}
\usepackage{ifthen}
\usepackage{cite}
\usepackage[cmex10]{amsmath} 

\usepackage{amssymb,amsfonts, amsthm}
\usepackage{textcomp}
\usepackage{xcolor}
\usepackage{nicefrac}
\usepackage{mathtools}

\newtheorem{theorem}{Theorem}
\newtheorem*{theorem*}{Theorem}
\newtheorem{proposition}{Proposition}
\newtheorem*{proposition*}{Proposition}
\newtheorem{lemma}{Lemma}

\newcommand{\bu}{{\mathbf{u}}}
\newcommand{\bU}{{\mathbf{U}}}
\newcommand{\bv}{{\mathbf{v}}}
\newcommand{\bV}{{\mathbf{V}}}

\newcommand{\bX}{{\mathbf{X}}}
\newcommand{\bx}{{\mathbf{x}}}
\newcommand{\by}{{\mathbf{y}}}
\newcommand{\bY}{{\mathbf{Y}}}
\newcommand{\bZ}{{\mathbf{Z}}}

\newcommand{\sT}{{\rm T}}
\newcommand\smallO{
	\mathchoice
	{{\scriptstyle\mathcal{O}}}
	{{\scriptstyle\mathcal{O}}}
	{{\scriptscriptstyle\mathcal{O}}}
	{\scalebox{.7}{$\scriptscriptstyle\mathcal{O}$}}
}
\newcommand{\R}{\mathbb{R}}
\newcommand{\cL}{\mathcal{L}}
\newcommand{\cN}{\mathcal{N}}
\newcommand{\cH}{\mathcal{H}}
\newcommand{\cZ}{\mathcal{Z}}
\newcommand{\E}{\mathbb{E}}
\newcommand{\Var}{{\mathbb{V}\mathrm{ar}}}
\newcommand{\simiid}[1][0pt]{\mathrel{\raisebox{#1}{$\sim$}}}
\newcommand{\iid}{\overset{\text{\tiny i.i.d.}}{\simiid[-2pt]}}

\begin{document}
\title{High-dimensional rank-one nonsymmetric matrix decomposition: the spherical case}

\author{
	\IEEEauthorblockN{
		Cl\'{e}ment Luneau\IEEEauthorrefmark{1}\IEEEauthorrefmark{3},
		Nicolas Macris\IEEEauthorrefmark{1}
		and Jean Barbier\IEEEauthorrefmark{2}
	}
	\IEEEauthorblockA{
		\IEEEauthorrefmark{1}Laboratoire de Th\'{e}orie des Communications,
		Ecole Polytechnique F\'{e}d\'{e}rale de Lausanne, Switzerland.
	}
	\IEEEauthorblockA{\IEEEauthorrefmark{2} The Abdus Salam International Center for Theoretical Physics, Trieste, Italy.}
	\IEEEauthorblockA{\IEEEauthorrefmark{3}Email: clement.luneau@epfl.ch}
}

\maketitle

\begin{abstract}
	We consider the problem of estimating a rank-one nonsymmetric matrix under additive white Gaussian noise.
	The matrix to estimate can be written as the outer product of two vectors and we look at the special case in which both vectors are uniformly distributed on spheres.
	We prove a replica-symmetric formula for the average mutual information between these vectors and the observations in the high-dimensional regime.
	This goes beyond previous results which considered vectors with independent and identically distributed elements. The method used can be extended to rank-one tensor problems.
\end{abstract}


\section{Introduction}\label{section:introduction}
Tensor decomposition, which originated with Hitchcock in 1927 \cite{Hitchcock1927}, has found many applications in signal processing, graph analysis, data mining and machine learning in the past two decades \cite{Sidiropoulos2017, cichoki2015, kolda2009}.
While tensor decomposition was originally developed in a deterministic and algebraic context, it is of interest for these applications to develop a statistical approach \cite{RicMon_2014}.
Some important questions in this setting are, for example, \textit{under which conditions and how can we recover a low-rank tensor -- the signal of interest -- from noisy observations of it?}
This work focuses on answering -- at least in part -- these questions in the most elementary, but yet rich, setting of a nonsymmetric rank-one matrix signal buried within noise.
Namely, we observe under additive white Gaussian noise (AWGN) a $n_u \times n_v$ rank-one matrix $\bU \bV^{\sT}$ where $\bU$ and $\bV$ are random vectors that we wish to recover as well as possible.
This problem, and its symmetric version, have generated important results in the past ten years \cite{Feral2007,BenaychGeorge2011}.

Our approach is in the continuity of a line of research establishing low-dimensional variational formulas for the average mutual information between a signal of interest and noisy observations in the high-dimensional regime \cite{Lesieur2017,Lelarge2019,Miolane2017nonsymmetric,Barbier2017phase}.
Such formulas are valuable because they link the mutual information of a high-dimensional channel whose outputs are coupled to those of simple decoupled scalar channels.
One can then determine, by solving a low-dimensional variational problem, phase transitions as well as performance measures related to the minimum mean square error (MMSE).
One also obtains important insights on the performance of (message passing) algorithms designed to estimate input signals.
In fact, the fixed points of the state evolution equations tracking the performance of the Approximate Message Passing algorithm in the high-dimensional regime can be identified among the critical points of the variational expression for the mutual information.

For the problem at hand, the variational formula -- that was predicted using the replica trick from statistical physics -- has already been proven rigorously \textit{when $\bU$ and $\bV$ have independent and identically distributed (i.i.d.) entries}\cite{Lelarge2019,Miolane2017nonsymmetric}.
These results were extended beyond the matrix case to rank-one nonsymmetric tensor decomposition \cite{BarbierMacrisAllerton2017,GanguliNIPS2018}.
The replica prediction has also been shown to be true for low-rank symmetric tensor decomposition \cite{Lelarge2019,luneauITW2019}.

A natural follow-up interrogation is what happens when either $\bU$ or $\bV$ doesn't have independent entries anymore. Can the average mutual information in the high-dimensional regime still be given by a \textit{simple}, \textit{low-dimensional}, variational formula?
In this work, we study the simple case in which both \textit{$\bU$ and $\bV$ are uniformly distributed on spheres} (whose radii scale like $\sqrt{n_u}$ and $\sqrt{n_v}$, respectively) and give a rigorous and positive answer to the question above. To the best of our knowledge fully rigorous results on this issue are scarce.
Recently, \cite{Aubin_NIPS2019} analyzed (under natural assumptions) another situation in which $\bU$ and $\bV$ are generated by a generalized linear model.

In Section~\ref{section:problem_setting} we present the problem and our main results.
In Section~\ref{section:proof_adaptive_interpolation} we give the reader an outline of the proof of the variational formula for the average mutual information.
We conclude in Section~\ref{section:conclusion} with a discussion of the relation between the present problem and the classical spherical spin-glass model of statistical mechanics.

\section{Problem setting and main results}\label{section:problem_setting}
Let $\bU \in \mathbb{R}^{n_u}$ and $\bV \in \R^{n_v}$ be uniformly distributed on the spheres of radii $\sqrt{\rho_u n_u}$ and $\sqrt{\rho_v n_v}$, respectively, with $\rho_u$ and $\rho_v$ positive real numbers.
We denote $P_u$ and $P_v$ their respective probability distributions.
The matrix factorization problem is the task of inferring both vectors $\bU$ and $\bV$ from a noisy observation of the scaled rank-one matrix $\bU \bV^{\sT}$. More precisely, we observe the matrix $\bY \in \mathbb{R}^{n_u \times n_v}$ whose entries satisfy $\forall (i,j) \in \{1,\dots,n_u\} \times \{1,\dots,n_v\}$:
\begin{equation}\label{def:Y}
Y_{ij} = \sqrt{\frac{\lambda}{n}}\,U_iV_j + Z_{ij}\;.
\end{equation}
Here, the matrix $\bZ \in \mathbb{R}^{n_u \times n_v}$ has i.i.d.\ elements with respect to (w.r.t.) the standard normal distribution $\cN(0, 1)$,
the positive real number $\lambda$ plays the role of a signal-to-noise ratio (SNR), and
the positive integer $n$ scales like $n_u$ and $n_v$, i.e., there exist positive real numbers $\alpha_u$ and $\alpha_v$ such that:
\begin{equation}\label{scaling_nu_and_nv}
\lim_{n\to + \infty} \frac{n_u}{n} = \alpha_u \;,\;
\lim_{n\to + \infty} \frac{n_v}{n} = \alpha_v \;.
\end{equation}
The normalization $1/\sqrt{n}$ in \eqref{def:Y} with the scaling \eqref{scaling_nu_and_nv} makes the estimation problem nontrivial.
Finally, we define the vector of hyperparameters for this problem: $\Theta \triangleq [\lambda \;\, \alpha_u \;\, \alpha_v \;\, \rho_u \;\, \rho_v]$.

\subsection{Variational formula for the average mutual information}
A central role is played by a simple linear model with AWGN and its average mutual information.
\begin{lemma}\label{lemma:free_entropy_spherical_vector}
Let $\bX$ be a $n$-dimensional random vector uniformly distributed on the sphere of radius $\sqrt{n}$ that is observed at the output of the following noisy linear channel:
\begin{equation}\label{spherical_vector_estimation}
\widetilde{\bY} = \sqrt{m}\,\bX + \widetilde{\bZ}
\end{equation}
where $\widetilde{Z}_i \iid \cN(0,1)$ for $i=1,\dots, n$ and $m > 0$ plays the role of a SNR.
The average mutual information between $\bX$ and $\widetilde{\bY}$ converges in the high-dimensional limit and:
\begin{equation}\label{avg_mutual_information_noisy_spherical_vector}
\lim_{n \to +\infty} \frac{I(\bX \,; \widetilde{\bY})}{n}
= \frac{\ln(1+m)}{2} \;.
\end{equation}
\end{lemma}

Note that the limit is equal to the average mutual information between $\bX$ and $\widetilde{\bY}$ where this time the entries of the signal $\bX$ are i.i.d.\ with respect to $\cN(0, 1)$.
It is well-known that such vector $\bX$ is approximately uniformly distributed on the sphere of radius $\sqrt{n}$ in high-dimension \cite[Section 3.3.3]{vershynin_2018}.
We now state our main theorem:
\begin{theorem}\label{theorem:limit_mutual_information}
Define the following \textit{potential function}:
\begin{equation}
i_{\scriptscriptstyle \Theta}(m_u,m_v)\triangleq \frac{\lambda \alpha_u \alpha_v}{2}(\rho_u-m_u)(\rho_v-m_v)
+ \alpha_u \frac{\ln(1+\lambda\alpha_v \rho_u m_v)}{2}
+ \alpha_v \frac{\ln(1 + \lambda \alpha_u \rho_v m_u)}{2}\:.\label{def:potential_mutual_information}
\end{equation}
In the high-dimensional limit, the average mutual information between $(\bU, \bV)$ and $\bY$ defined in \eqref{def:Y} satisfies:
\begin{equation}\label{limit_mutual_information}
\lim_{n\to +\infty} \frac{I(\bU,\bV\,;\,\bY)}{n}
= \adjustlimits{\inf}_{m_u \in [0,\rho_u]}{\sup}_{m_v \in [0,\rho_v]}\, i_{\scriptscriptstyle \Theta}(m_u,m_v) \;.
\end{equation}
\end{theorem}

The proof of Theorem~\ref{theorem:limit_mutual_information} is given in Section~\ref{section:proof_adaptive_interpolation}.
Note that the last two summands in \eqref{def:potential_mutual_information} are the asymptotic average mutual informations of two decoupled channels:
\begin{align*}
\lim_{n \to +\infty} I(\bU \,; \sqrt{\lambda\alpha_v m_v}\bU + \widetilde{\bZ})/n &= \alpha_u \ln(1+\lambda\alpha_v \rho_u m_v)/2\,;\\
\lim_{n \to +\infty} I(\bV \,; \sqrt{\lambda\alpha_u m_u}\bV + \overline{\bZ})/n &= \alpha_v \ln(1 + \lambda \alpha_u \rho_v m_u)/2\,;
\end{align*}
where $\widetilde{\bZ}$, $\overline{\bZ}$ are independent AWGN.

We remark that the limit of $\nicefrac{I(\bU,\bV;\bY)}{n}$ is the same if both $\bU$ and $\bV$ have i.i.d.\ standard Gaussian components (see \cite{BarbierMacrisAllerton2017,Miolane2017nonsymmetric}).
The equivalence of the spherical and Gaussian cases is not 
an obvious fact when it comes to make a precise argument.
We discuss this point further in Section~\ref{section:conclusion}.

\subsection{Minimum mean square error}
It is well-known that the mean square error of an estimator of $\bU \bV^T$ that is a function of $\bY$ \textit{only} is minimized by the posterior mean $\E[\bU \bV^T \vert \bY]$.
We denote by $\mathrm{MMSE}_\lambda(\bU\bV^T \vert \bY)$ the minimum mean square error
$\nicefrac{\E\,\Vert \bU\bV^T - \E[\bU \bV^T \vert \bY] \Vert^2 }{n_u n_v}$ (it depends on $\lambda$ through the observations $\bY$).
Combining  Theorem~\ref{theorem:limit_mutual_information} with the {I-MMSE} relation (see \cite{Guo2005_IMMSE})
\begin{equation*}
\frac{\partial}{\partial \lambda}\bigg(\frac{I(\bU, \bV ; \bY)}{n}\bigg)
= \frac{n_u}{n}\frac{n_v}{n} \frac{\mathrm{MMSE}_\lambda(\bU\bV^T \vert \bY)}{2}
\end{equation*}
yields the next theorem. Its proof is given in Appendix~\ref{app:proof_mmse}.
\begin{theorem}\label{theorem:MMSE}
Let $\lambda_{\mathrm{IT}} \triangleq \frac{1}{\rho_u\rho_v \sqrt{\alpha_u\alpha_v}}$.
For all $\lambda \in (0,+\infty)$, there is a unique solution to the extremization over $(m_u,m_v)$ on the right-hand side of \eqref{limit_mutual_information} given by:
\begin{equation*}
\big(m_u^*(\lambda),m_v^*(\lambda) \big)
=
\begin{cases}
\qquad\qquad\qquad\;\;\: (0\,,0)&\text{if}\;\: 0 < \lambda \leq \lambda_{\mathrm{IT}}\\
\Big(\frac{\lambda^2 \alpha_u \alpha_v \rho_v^2 \rho_u^2 - 1}{\lambda \alpha_u \rho_v (1 + \lambda \alpha_v \rho_v \rho_u)}
,\!\frac{\lambda^2 \alpha_u \alpha_v \rho_v^2 \rho_u^2 - 1}{\lambda \alpha_v \rho_u (1 + \lambda \alpha_u \rho_v \rho_u)}\Big)&\text{if}\;\: \lambda > \lambda_{\mathrm{IT}}
\end{cases}.
\end{equation*}
Then, $\mathrm{MMSE}_\lambda(\bU\bV^T \vert \bY)$ satisfies:
	\begin{equation}
	\lim_{n \to +\infty} \mathrm{MMSE}_\lambda(\bU\bV^T \vert \bY) = \rho_u \rho_v - m_u^*(\lambda) m_v^*(\lambda)\;.
	\end{equation}
Hence, the asymptotic MMSE is less than $\rho_u \rho_v$ if, and only if, $\lambda > \lambda_{\mathrm{IT}}$.
\end{theorem}

Theorems \ref{theorem:limit_mutual_information} and \ref{theorem:MMSE}  provide important insight on the inference problem.
Nonanalytic points of \eqref{limit_mutual_information} correspond to the location of \textit{phase transitions} where the MMSE changes behavior.
In the present problem, we find by an explicit analysis a unique continuous phase transition point $\lambda_{\rm{IT}}$.
The mutual information is continuously differentiable for all $\lambda>0$ and its second derivative has a jump at $\lambda_{\mathrm{IT}}$.
Correspondingly, the MMSE is continuous with a jump in its first derivative at $\lambda_{\mathrm{IT}}$.
More precisely, the MMSE is $\rho_u \rho_v$ for $\lambda \leq \lambda_{\mathrm{IT}}$ and it continuously departs from $\rho_u \rho_v$ once $\lambda$ becomes greater than $\lambda_{\mathrm{IT}}$.
Thus, $\lambda_{\rm IT}$ is the lowest SNR for which an estimate of the matrix $\bU \bV^T$ is \textit{information-theoretically} possible.
The general phenomenological picture has been uncovered in a number of situations (including richer ones) by direct analysis of the replica formula for the asymptotic mutual information.
We refer to \cite{Lesieur2017} for more details.


\section{Proof of Theorem~\ref{theorem:limit_mutual_information}}\label{section:proof_adaptive_interpolation}
We only present the main ideas and steps of the proof.
We will refer to the appendices, which contain all the technicalities of the proof, when needed.
The proof is based on the adaptive interpolation method introduced in \cite{BarbierMacrisAdaptiveInterpolation, BarbierMacris2019}.
The main difference with the canonical interpolation method developed by Guerra and Toninelli in the context of spin glasses \cite{GuerraToninelli2002, Guerra2003} is the increased flexibility in choosing the path followed by the interpolation between its two extremes.
By choosing two different interpolation paths, we will bound the asymptotic average mutual information from above and below by the same variational formula.
For the proof we assume that $\lambda = 1$.
This is without loss of generality as we can always reduce to this case by rescaling $\rho_u$ to $\lambda \rho_u$.
\subsection{Adaptive path interpolation}
We introduce a ``time'' parameter $t \in [0,1]$.
The adaptive interpolation interpolates from the original channel \eqref{def:Y} at $t=0$ to two independent channels similar to \eqref{spherical_vector_estimation} at $t=1$ (one for $\bU$ and one for $\bV$).
In between, we follow an interpolation path $R(\cdot,\epsilon) = (R_u(\cdot,\epsilon), R_v(\cdot,\epsilon))$ where $R_u(\cdot,\epsilon)$ and $R_v(\cdot,\epsilon)$ are continuously differentiable functions from $[0,1]$ to $[0,+\infty)$ parametrized by a ``small perturbation'' $\epsilon =(\epsilon_u, \epsilon_v) \in [0,+\infty)^2$ and such that $R(0,\epsilon)=\epsilon$.
More precisely, for $t \in [0,1]$, we observe:
\begin{align}\label{interpolation_model}
\begin{cases}
\bY^{(t)} &= \sqrt{\frac{1-t}{n}}\,\bU \, \bV^{\sT} + \bZ\quad\;;\\
\widetilde{\bY}^{(t,\epsilon)} &= \sqrt{\alpha_v R_v(t,\epsilon)}\, \bU + \widetilde{\bZ} \;\;\:;\\
\overline{\bY}^{(t,\epsilon)} &= \sqrt{\alpha_u R_u(t,\epsilon)}\, \bV + \overline{\bZ}\;\,;
\end{cases}
\end{align}
where $\bU \sim P_u$, $\bV \sim P_v$ and all of the noises $\bZ \in \R^{n_u \times n_v}$, $\widetilde{\bZ} \in \R^{n_u}$, $\overline{\bZ} \in \R^{n_v}$ have i.i.d.\ entries with respect to $\cN(0,1)$.
Applying Bayes' rule, we obtain the posterior distribution of $(\bU, \bV)$ given $(\bY^{(t)},\widetilde{\bY}^{(t,\epsilon)}, \overline{\bY}^{(t,\epsilon)})$:
\begin{equation}\label{posterior_distribution_U_V}
dP(\bu, \bv \vert \bY^{(t)},\widetilde{\bY}^{(t,\epsilon)}, \overline{\bY}^{(t,\epsilon)})
\triangleq \frac{dP_{u}(\bu)dP_v(\bv) \, e^{-\cH_{t,\epsilon}(\bu, \bv ; \bY^{(t)},\widetilde{\bY}^{(t,\epsilon)}, \overline{\bY}^{(t,\epsilon)})}}{\cZ_{t,\epsilon}(\bY^{(t)},\widetilde{\bY}^{(t,\epsilon)}, \overline{\bY}^{(t,\epsilon)})} \;,
\end{equation}
where we introduced the interpolating Hamiltonian
\begin{align}\label{interpolating_hamiltonian}
\cH_{t,\epsilon}(\bu, \bv ; \bY^{(t)}, \widetilde{\bY}^{(t,\epsilon)}, \overline{\bY}^{(t,\epsilon)})
&\triangleq
\sum_{i=1}^{n_u}\sum_{j=1}^{n_v}
\frac{1-t}{2n} u_i^2 v_j^2 - \sqrt{\frac{1-t}{n}}\, u_i v_j Y_{ij}^{(t)}
+
\sum_{i=1}^{n_u} \frac{\alpha_v R_v(t,\epsilon)}{2} u_i^2 - \sqrt{\alpha_v R_v(t,\epsilon)}\, u_i \widetilde{Y}_i^{(t,\epsilon)} \nonumber\\
&\qquad\qquad\qquad\qquad\qquad\qquad\qquad\qquad\quad\:
+\sum_{j=1}^{n_v} \frac{\alpha_u R_u(t,\epsilon)}{2} v_j^2 - \sqrt{\alpha_u R_u(t,\epsilon)}\, v_j \overline{Y}_j^{(t,\epsilon)} \:,
\end{align}
and $\cZ_{t,\epsilon}(\bY^{(t)},\widetilde{\bY}^{(t,\epsilon)}, \overline{\bY}^{(t,\epsilon)})$ properly normalizes the posterior.
Note that \eqref{interpolating_hamiltonian} could be simplified using the spherical constraints but this general form is convenient for the analysis.
We denote an expectation with respect to the posterior distribution \eqref{posterior_distribution_U_V} using the angular brackets $\langle - \rangle_{t,\epsilon}$, i.e., 
$\langle g(\bu, \bv) \rangle_{t,\epsilon} = \int \!g(\bu, \bv)\,dP(\bu, \bv \vert \bY^{(t)},\widetilde{\bY}^{(t,\epsilon)}, \overline{\bY}^{(t,\epsilon)})$.
The interpolating average free entropy defined as
\begin{equation}\label{interpolating_free_entropy}
f_n(t,\epsilon) \triangleq \frac1n \E \ln \cZ_{t,\epsilon}(\bY^{(t)},\widetilde{\bY}^{(t,\epsilon)}, \overline{\bY}^{(t,\epsilon)})
\end{equation}
is intimately linked to the average mutual information.
In particular, $f_n \triangleq f_n(0,0) = \nicefrac{n_u n_v \rho_u \rho_v}{2n^2}-\nicefrac{I(\bU, \bV;\bY)}{n}$.
Hence, Theorem~\ref{theorem:limit_mutual_information} is equivalent to:
\begin{equation}\label{limit_free_entropy}
\lim_{n\to \infty} f_n\\
= \adjustlimits{\sup}_{m_u}{\inf}_{m_v}\: \phi_{\scriptscriptstyle \Theta}(m_u,m_v)\;,
\end{equation}
where the potential $\phi_{\scriptscriptstyle \Theta}$ is defined using $\varphi(m) \triangleq \frac{m - \ln(1+m)}{2}$:
\begin{equation*}
\phi_{\scriptscriptstyle \Theta}(m_u,\!m_v) \!\triangleq\! 
\alpha_u \varphi(\alpha_v \rho_u m_v) + \alpha_v \varphi(\alpha_u \rho_v m_u)
-\frac{\alpha_u \alpha_v m_u m_v}{2}\;.
\end{equation*}
Looking at how $f_n(t,\epsilon)$ varies from $t=0$ to $t=1$ yields the following important \textit{sum-rule} that we will later evaluate for different interpolation paths.
\begin{proposition}\label{prop:sum_rule}
Define the scalar overlaps
$Q_u \triangleq \frac{1}{n_u} \sum_i u_i U_i$ and
$Q_v \triangleq \frac{1}{n_v} \sum_i v_i V_i$.
Denote $R'_u(\cdot,\epsilon)$ and $R'_v(\cdot,\epsilon)$ the derivative of $R_u(\cdot,\epsilon)$ and $R_v(\cdot,\epsilon)$, respectively.
Assume that both $R'_u(t,\epsilon)$ and $R'_v(t,\epsilon)$ are uniformly bounded in $(t,\epsilon)$ belonging to $[0,1] \times [0,+\infty)^2$.
Then:
\begin{align}
f_n = \mathcal{O}(\Vert \epsilon \Vert) + \smallO_n(1) + \alpha_u \varphi(\alpha_v \rho_u R_v(1,\epsilon))
+ \alpha_v \varphi(\alpha_u \rho_v R_u(1,\epsilon))
&-\frac{\alpha_u \alpha_v}{2} \!\! \int_{0}^{1} \!\!\! dt\, R'_u(t,\epsilon) R'_v(t,\epsilon)\nonumber\\
&+ \frac{\alpha_u \alpha_v}{2} \!\! \int_{0}^{1} \!\!\! dt\,
\E\big\langle (Q_u \! - \! R'_u(t,\epsilon)) (Q_v \! - \! R'_v(t,\epsilon))\big\rangle_{t,\epsilon}
\label{sum_rule}
\end{align}
where $\smallO_n(1)$ is a quantity that vanishes uniformly in $\epsilon$ as $n$ gets large,
and $\mathcal{O}(\Vert \epsilon \Vert)$ is a quantity whose absolute value is upper bounded by $C \Vert \epsilon \Vert$ for some constant $C$ independent of both $n$ and $\epsilon$.
\end{proposition}
\begin{IEEEproof}
Evaluating \eqref{interpolating_free_entropy} at both extremes of the interpolation yields
$f_n(0,\epsilon) = f_n(0,0) + \mathcal{O}(\Vert \epsilon \Vert) = f_n + \mathcal{O}(\Vert \epsilon \Vert)$ and 
$f_n(1,\epsilon) = \alpha_u \varphi(\alpha_v \rho_u R_v(1,\epsilon)) + \alpha_v \varphi(\alpha_u \rho_v R_u(1,\epsilon)) + \smallO_n(1)$ where $\smallO_n(1)$, $\mathcal{O}(\Vert \epsilon \Vert)$ are quantities satisfying the properties given in the proposition.
We obtain the sum-rule \eqref{sum_rule} by combining the later with the fundamental theorem of calculus
$f_n(0,\epsilon)=f_n(1,\epsilon)-\int_{0}^{1}f_n^{\prime}(t,\epsilon) dt$,
where $f_n^{\prime}(\cdot,\epsilon)$ is the derivative of $f_n(\cdot,\epsilon)$.
All the technical details, including the computation of $f_n^{\prime}(\cdot,\epsilon)$, is given in Appendix~\ref{app:sum_rule}.
\end{IEEEproof}
\subsection{Interpolation paths as solutions to ODEs}
To prove Theorem~\ref{theorem:limit_mutual_information}, we will lower bound $\liminf_n f_n$ and upper bound $\limsup_n f_n$ by the same quantity $\sup_{m_u}\inf_{m_v}\phi_{\scriptscriptstyle \Theta}(m_u,m_v)$.
To do so we will plug two different choices for $R(\cdot,\epsilon)$ in the sum-rule \eqref{sum_rule}.
In both cases, $R(\cdot,\epsilon)$ will be the solution of a second-order ordinary differential equation (ODE). We now describe these ODEs before diving further into the proofs of the matching bounds.

For $t \in [0,1]$ and $R = (R_u, R_v) \in [0, +\infty)^2$, consider the problem of estimating $(\bU, \bV)$ from the observations:
\begin{align}\label{interpolation_model_R}
\begin{cases}
\bY^{(t)} &= \sqrt{\frac{1-t}{n}}\,\bU \, \bV^{\sT} + \bZ\;;\\
\widetilde{\bY}^{(t,R_v)} &= \sqrt{\alpha_v R_v}\, \bU + \widetilde{\bZ} \;\;\:;\\
\overline{\bY}^{(t,R_u)} &= \sqrt{\alpha_u R_u}\, \bV + \overline{\bZ}\;\;;
\end{cases}
\end{align}
where $\bU \sim P_u$, $\bV \sim P_v$ and all of the noises $\bZ \in \R^{n_u \times n_v}$, $\widetilde{\bZ} \in \R^{n_u}$, $\overline{\bZ} \in \R^{n_v}$ have i.i.d.\ entries with respect to $\cN(0,1)$.
The posterior distribution of $(\bU, \bV)$ given $(\bY^{(t)},\widetilde{\bY}^{(t,R_v)}, \overline{\bY}^{(t,R_u)})$ is (up to the normalization factor):
\begin{equation}\label{posterior_H_t_R}
dP(\bu, \bv \vert \bY^{(t)},\widetilde{\bY}^{(t,R_v)}, \overline{\bY}^{(t,R_u)})
\propto
dP_{u}(\bu)dP_v(\bv) \, e^{-\cH_{t,R}(\bu, \bv ; \bY^{(t)},\widetilde{\bY}^{(t,R_v)}, \overline{\bY}^{(t,R_u)})} \;,
\end{equation}
where $\cH_{t,R}$ denotes the associated interpolating Hamiltonian:
\begin{align*}
\cH_{t,R}(\bu, \bv ; \bY^{(t)}, \widetilde{\bY}^{(t,R_v)}, \overline{\bY}^{(t,R_u)})
\triangleq
\sum_{i=1}^{n_u}\sum_{j=1}^{n_v}
\frac{1-t}{2n} u_i^2 v_j^2 - \sqrt{\frac{1-t}{n}}\, u_i v_j Y_{ij}^{(t)}
&+ \sum_{i=1}^{n_u} \frac{\alpha_v R_v}{2} u_i^2 - \sqrt{\alpha_v R_v}\, u_i \widetilde{Y}_i^{(t,R_v)}\\
&+ \sum_{j=1}^{n_v} \frac{\alpha_u R_u}{2} v_j^2 - \sqrt{\alpha_u R_u}\, v_j \overline{Y}_j^{(t,R_u)} \;.
\end{align*}
The angular brackets $\langle - \rangle_{t,R}$ will denote the expectation w.r.t.\ the posterior \eqref{posterior_H_t_R}.
Let $m_u \in [0,\rho_u]$, $F_v(t, R) \triangleq \E \langle Q_v \rangle_{t,R}$ and $F_u(t, R) \triangleq 2 \rho_u \varphi^{\prime}(\alpha_v \rho_u \E \langle Q_v \rangle_{t,R})$. We will consider the two following second-order ODEs with initial value $\epsilon \in [0,+\infty)^2$:
\begin{align}
&g' = \;\big(m_u , F_v(t,g)\big)\quad\;\;\, , \: g(0)=\epsilon \;; \label{def:ode_lowerbound}\\
&g' = \big(F_u(t,g), F_v(t, g)\big) \, , \: g(0)=\epsilon\;.\label{def:ode_upperbound}
\end{align}
The next proposition sums up useful properties on the solutions of these two ODEs, i.e., our two kinds of interpolation paths.
The proof is given in Appendix~\ref{app:properties_ode}.
\begin{proposition}\label{prop:ode}
	For all $\epsilon \in [0,+\infty)^2$, there exists a unique global solution $R(\cdot,\epsilon): [0,1] \to [0,+\infty)^2 $ to \eqref{def:ode_lowerbound}.
	This solution is continuously differentiable and its derivative $R'(\cdot,\epsilon)$ satisfies $R'([0,1],\epsilon) \subseteq [0, \rho_u] \times [0, \rho_v]$.
	Besides, for all $t \in [0,1]$, $R(t,\cdot)$ is a $\mathcal{C}^1$-diffeomorphism from $[0, +\infty)^2$ into its image whose Jacobian determinant is greater than, or equal to, one:
	\begin{equation}\label{jacobian_greater_one_ode_upperbound}
	\forall \,\epsilon \in [0, +\infty)^2: \det J_{R(t,\cdot)}(\epsilon) \geq 1 \:,
	\end{equation}
	where $J_{R(t,\cdot)}$ denotes the Jacobian matrix of $R(t,\cdot)$.\\
	Finally, the same statement holds if we consider \eqref{def:ode_upperbound} instead.
\end{proposition}
\subsection{Proof of the lower bound on $\liminf_n f_n$}
Let $m_u \in [0, \rho_u]$ and $\epsilon = (\epsilon_u, \epsilon_v) \in (0,+\infty)^2$.
We choose as interpolation path the unique solution $R(\cdot, \epsilon)$ to the ODE \eqref{def:ode_lowerbound}.
Then $R'_u(t,\epsilon) = m_u$ and $R'_v(t,\epsilon) = \E \langle Q_v \rangle_{t, \epsilon}$.
Plugging this choice in the sum-rule of Proposition~\ref{prop:sum_rule}, and making use of
$\varphi(\alpha_v \rho_u R_v(1,\epsilon)) = \mathcal{O}(\epsilon_v) + \varphi\big(\alpha_v \rho_u \int_0^1 dt \, R'_v(1,\epsilon)\big)$ and $\varphi(\alpha_u \rho_v (\epsilon_u + m_u)) = \mathcal{O}(\epsilon_u) + \varphi(\alpha_u \rho_v m_u)$, yields:
\begin{align}
f_n
&= \mathcal{O}(\Vert \epsilon \Vert) + \smallO_n(1) + \phi_{\scriptscriptstyle \Theta}\Big(m_u,\int_0^1 dt \, R'_v(1,\epsilon)\Big)
+ \frac{\alpha_u \alpha_v}{2} \int_{0}^{1} dt\,
\E\,\big\langle Q_u \big(Q_v - \E \langle Q_v \rangle_{t, \epsilon}\big)\big\rangle_{t,\epsilon}\nonumber\\
&\geq \mathcal{O}(\Vert \epsilon \Vert) + \smallO_n(1) + \; \inf_{m_v \in [0, \rho_v]} \phi_{\scriptscriptstyle \Theta}(m_u,m_v)\;
+ \frac{\alpha_u \alpha_v}{2} \int_{0}^{1} dt\,
\E\,\big\langle Q_u \big(Q_v - \E \langle Q_v \rangle_{t, \epsilon}\big)\big\rangle_{t,\epsilon}\;. \label{sumrule_lowerbound}
\end{align}
The lower bound is because $\int_0^1 dt \, R'_v(1,\epsilon) \in [0, \rho_v]$ (see Proposition~\ref{prop:sum_rule}).
If the overlap $Q_v$ concentrates on its expectation then the remainder $\int_{0}^{1} dt\,
\E\langle Q_u (Q_v - \E \langle Q_v \rangle_{t, \epsilon})\rangle_{t,\epsilon}$ in the lower bound \eqref{sumrule_lowerbound} vanishes and we can end the proof.
However, proving such concentration is only possible after integrating on a well-chosen set of ``perturbation'' $\epsilon$.
This integration over $\epsilon$ \textit{smoothens} the phase transitions that might appear for particular choices of $\epsilon$ when $n$ goes to infinity.
From now on, $\epsilon \in \mathcal{S}_n \triangleq [s_n, 2s_n]^2$ where $s_n \triangleq n^{-\eta}$, $\eta > 0$.
Integrating w.r.t. \ $\epsilon$ on both sides of \eqref{sumrule_lowerbound} yields $\big(f_n = \int_{[s_n,2s_n]^2} \nicefrac{f_nd\epsilon}{s_n^2}\big)$:
\begin{equation}
f_n
\geq
\smallO_n(1) + \inf_{m_v \in [0, \rho_v]} \phi_{\scriptscriptstyle \Theta}(m_u,m_v) + \frac{\alpha_u \alpha_v}{2}\mathcal{R}\;,\label{sumrule_lowerbound_integration_epsilon}
\end{equation}
where $\mathcal{R} \triangleq \int_{0}^{1} dt \int_{\mathcal{S}_n} \frac{d\epsilon}{s_n^2}
\E\,\langle Q_u (Q_v - \E \langle Q_v \rangle_{t, \epsilon_v})\rangle_{t,\epsilon}$.
By Jensen's inequality and $\vert Q_u \vert \leq \nicefrac{\Vert \bU \Vert \Vert \bu \Vert}{n_u} = \rho_u$, $\mathcal{R}$ satisfies:
\begin{equation}\label{upperbound_integrated_remainder}
\vert \mathcal{R} \vert \leq
\rho_u\int_{0}^{1} dt\,
\sqrt{\int_{\mathcal{S}_n} \frac{d\epsilon}{s_n^2} \, \E\,\big\langle (Q_v - \E \langle Q_v \rangle_{t, \epsilon})^2\big\rangle_{t,\epsilon}} \;.
\end{equation}
The change of variables $\epsilon \to R = R(t,\epsilon)$ -- justified by $R(t,\cdot)$ being a $\mathcal{C}^1$-diffeomorphism from $[0,+\infty)$ to its image (see Proposition~\ref{prop:ode}) -- yields for all $t \in [0,1]$:
\begin{equation*}
\int_{\mathcal{S}_n} \frac{d\epsilon}{s_n^2} \E\,\big\langle (Q_v - \E \langle Q_v \rangle_{t, \epsilon_v})^2\big\rangle_{t,\epsilon}
\leq \int_{s_n}^{2s_n + \rho_v}  \! \frac{dR_v}{s_n^2}  \int_{s_n}^{2s_n + \rho_u} \! dR_u \,\E\,\big\langle (Q_v - \E \langle Q_v \rangle_{t, R})^2\big\rangle_{t,R} \;. 
\end{equation*}
The inequality follows from the integrand being nonnegative, $\eqref{jacobian_greater_one_ode_upperbound}$ and $R(t,\mathcal{S}_n) \subseteq [s_n, 2s_n + \rho_u] \times [s_n, 2s_n + \rho_v]$.
We now apply Proposition~\ref{prop:concentration_overlap} -- an important result on the concentration of the overlap $Q_v$ that follows this proof -- with $M_u = 2 + \rho_u$, $M_v = 2 + \rho_v$, $a=s_n$, $b=2s_n + \rho_u$ and $\delta = s_n n^{\frac{2\eta-1}{3}}$ (we further assume $\eta < 1/2$).
Then, for $n$ large enough, there exists $M > 0$ such that $\forall t \in [0,1]$, $\forall R_v \in [s_n, 2s_n + \rho_v]$:
\begin{equation*}
\int_{s_n}^{2s_n + \rho_u} dR_u \,\E\,\big\langle (Q_v - \E \langle Q_v \rangle_{t, R})^2\big\rangle_{t,R}
\leq M / n^{\frac{1-2\eta}{3}}\;.
\end{equation*}
Combining this inequality with the two previous ones, we get $\vert \mathcal{R} \vert
\leq M' n^{\frac{4\eta}{3} -\frac{1}{6}}$ where $M' \triangleq \rho_u \sqrt{(1 + \rho_v)M}$.
This upper bound on $\vert \mathcal{R} \vert$ vanishes for $n$ large as long as $\eta$ is less than $\nicefrac{1}{8}$.
Passing to the limit inferior on both sides of \eqref{sumrule_lowerbound_integration_epsilon} thus yields
$\liminf_{n \to +\infty} f_n \geq
\inf_{m_v \in [0, \rho_v]} \phi_{\scriptscriptstyle \Theta}(m_u,m_v)$.
As this is true for all $m_u \in [0, \rho_u]$, we finally obtain:
$$
\liminf_{n \to +\infty} f_n \geq
\sup_{m_u \in [0,\rho_u]}\inf_{m_v \in [0, \rho_v]} \phi_{\scriptscriptstyle \Theta}(m_u,m_v)\;.
$$
\hfill\IEEEQEDclosed
\subsection{Concentration of the overlap $Q_v$}
We rely on the following concentration result to prove the matching bounds. It is clear that a similar result holds for $Q_u$.
\begin{proposition}\label{prop:concentration_overlap}
	Let $M_u, M_v >0$.
	For $n$ large enough, there exists a constant $M$ such that $\forall (a,b) \in (0,M_u)^2$ with ${a < \min\{1,b\}}$, $\forall \delta \in (0,a)$, $\forall R_v \in [0,M_v]$, $\forall t \in [0,1]$:
	\begin{equation*}
	\int_{a}^{b} \E\,\big\langle (Q_v -\E\,\langle Q_v \rangle_{t,R})^2\,\big\rangle_{t,R}\, dR_u
	\leq M\bigg(\frac{1}{\delta^2 n} - \frac{\ln(a)}{n} + \frac{\delta}{a-\delta}\bigg)\;.
	\end{equation*}
\end{proposition}
The proof, technical but not difficult, is given in Appendix~\ref{app:concentration_overlap}, and
follows the same step than similar concentration results on the overlaps of inference problems \cite{BarbierMacrisAdaptiveInterpolation, BarbierMacris2019, BarbierMacrisAllerton2017, Barbier2017phase}.
The differences with the proof in \cite{BarbierMacrisAllerton2017} are due to the entries of both $\bU$ and $\bV$ being not independent anymore.
It mainly impacts the proof that the free entropy $\nicefrac{\ln \cZ_{t,R}}{n}$\footnote{ $\cZ_{t,R} \equiv \cZ_{t,R}(\bY^{(t)},\widetilde{\bY}^{(t,R_v)}, \overline{\bY}^{(t,R_u)})$ is the normalization factor of the right-hand side of \eqref{posterior_H_t_R}.} concentrates on its average, which we need in our proof of the overlap concentration. We now use L\'{e}vy's lemma \cite[Corollary 5.4]{vershynin_2018} to show that $\nicefrac{\ln \cZ_{t,R}}{n}$ concentrates on its expectation with respect to $\bU, \bV$.
This requires verifying that $(\widetilde{\bU}, \widetilde{\bV}) \mapsto \nicefrac{\ln \cZ_{t,R}}{n}$ is Lipschitz continuous with respect to $\widetilde{\bU} \triangleq \nicefrac{\bU}{\sqrt{\rho_u n_u}}$ on the {$(n_u-1)$-sphere} and $\widetilde{\bV} \triangleq \nicefrac{\bV}{\sqrt{\rho_v n_v}}$ on the $(n_v-1)$-sphere.
The other difference is that in \cite[Lemma 3.1]{BarbierMacrisAllerton2017} the concentration result holds under the assumption that the prior of the i.i.d.\ entries of $\bV$ is compactly supported.
Here, knowing that the norm of $\bV$ scales like $\sqrt{n}$ is in fact enough to guarantee Proposition~\ref{prop:concentration_overlap}.

\subsection{Proof of the matching upper bound on $\limsup_n f_n$}
Let $\epsilon = (\epsilon_u, \epsilon_v) \in (0, +\infty)^2$.
We choose as interpolation path the unique solution $R(\cdot, \epsilon)$ to the ODE \eqref{def:ode_upperbound}.
Then, $R'_u(t,\epsilon) = 2 \rho_u \varphi'(\alpha_v \rho_u \E \langle Q_v \rangle_{t, \epsilon})$ and $R'_v(t,\epsilon) = \E \langle Q_v \rangle_{t, \epsilon}$.
Note that $\varphi(\alpha_v \rho_u \int_0^1 dt\, R'_v(t,\epsilon)) \leq \int_0^1  dt\, \varphi(\alpha_v \rho_u R'_v(t,\epsilon))$ as $\varphi$ is convex. So
$$
\varphi(\alpha_v \rho_u R_v(1,\epsilon)) \leq \mathcal{O}(\epsilon_v) + \int_0^1 dt\, \varphi(\alpha_v \rho_u R'_v(t,\epsilon))\;,
$$
and a similar inequality holds for $\varphi(\alpha_u \rho_v R_u(1,\epsilon))$.
Making use of these bounds after writing the sum-rule of Proposition~\ref{prop:sum_rule} for this particular interpolation path gives:
\begin{align}
f_n
\leq \mathcal{O}(\Vert \epsilon \Vert) + \smallO_n(1)
+ \int_0^1 dt \,\phi\big(R'_u(t,\epsilon), R'_v(t,\epsilon)\big)
+ \frac{\alpha_u \alpha_v}{2}  \int_{0}^{1}  dt\,
\E\,\big\langle(Q_u  - R'_u(t,\epsilon))(Q_v - \E \langle Q_v \rangle_{t, \epsilon})\big\rangle_{t,\epsilon} \;.\label{sumrule_upperbound}
\end{align}
Fix $t \in [0,1]$ and let ${h: m_v \in [0, \rho_v] \mapsto \phi(R'_u(t,\epsilon), m_v)}$.
As $R'_u(t,\epsilon) = 2 \rho_u \varphi'(\alpha_v \rho_u R'_v(t,\epsilon))$, we have $h'(R'_v(t,\epsilon)) = 0$
and the unique global minima of the strictly convex function $h$ is reached at $R'_v(t,\epsilon) \in [0,\rho_v]$, i.e.,
$$
\phi(R'_u(t,\epsilon), R'_v(t,\epsilon)) = \inf_{m_v \in [0,\rho_v]} \phi(R'_u(t,\epsilon), m_v)\;.
$$
Hence $\phi(R'_u(t,\epsilon), R'_v(t,\epsilon)) \leq \sup_{m_u} \inf_{m_v} \phi(m_u, m_v)$.
Plugging this upper bound back in \eqref{sumrule_upperbound} yields:
\begin{equation}
f_n
\leq \mathcal{O}(\Vert \epsilon \Vert) + \smallO_n(1) + {\adjustlimits \sup_{m_u \in [0, \rho_u]} \inf_{m_v \in [0, \rho_v]}}
\phi\big(m_u, m_v\big)
+ \frac{\alpha_u \alpha_v}{2}\! \int_{0}^{1} dt\,
\E\,\big\langle (Q_u  -  R'_u(t,\epsilon))(Q_v  - \E \langle Q_v \rangle_{t, \epsilon})\big\rangle_{t,\epsilon} \;.\label{upperbound_f_n_with_remainder}
\end{equation}
We get rid of the remainder exactly as in the proof of the lower bound.
After integrating \eqref{upperbound_f_n_with_remainder} over $\epsilon \in \mathcal{S}_n \triangleq [s_n, 2 s_n]^2$ $(s_n \triangleq n^{-\eta}$ with $\eta > 0)$, $f_n = \int_{\mathcal{S}_n} \nicefrac{f_n \, d\epsilon}{s_n^2}$ satisfies:
\begin{equation}\label{upperbound_fn_remainder_R}
f_n \leq \smallO_n(1) + \sup_{m_u} \inf_{m_v} \phi(m_u, m_v) + \nicefrac{\alpha_u \alpha_v \mathcal{R}}{2} \;,
\end{equation}
where $\mathcal{R}$ stands for the remainder:
\begin{equation*}
\mathcal{R} \triangleq \int_{0}^{1} dt \int_{\mathcal{S}_n} \frac{d\epsilon}{s_n^2}\,
\E\,\big\langle (Q_u - R'_u(t,\epsilon))(Q_v - \E \langle Q_v \rangle_{t, \epsilon})\big\rangle_{t,\epsilon} \;.
\end{equation*}
We can upper bound the absolute value of $\mathcal{R}$ by $C n^{\frac{4\eta}{3}-\frac{1}{6}}$ for some positive constant $C$ and $n$ large enough.
It is done exactly as in the proof of the lower bound on $\liminf_n f_n$:
$\vert Q_u - R'_u(t,\epsilon)\vert \leq 2\rho_u$ and the change of variables $\epsilon \to R = R(t,\epsilon)$ is justified by $R(t,\cdot)$ being a $\mathcal{C}^1$-diffeomorphism from $[0,+\infty)$ to its image (see Proposition~\ref{prop:ode}).
As long as $\eta$ is less than $\nicefrac{1}{8}$, the remainder vanishes when $n$ goes to infinity and passing to the limit superior on both sides of the inequality \eqref{upperbound_fn_remainder_R} yields the desired upper bound:
\begin{equation}
\limsup_{n \to +\infty} f_n \leq
\sup_{m_u \in [0, \rho_u]} \inf_{m_v \in [0, \rho_v]} \: \phi(m_u, m_v) \;.
\end{equation}
\hfill\IEEEQEDclosed

\section{Conclusion}\label{section:conclusion}
We conclude with a few comments on close connections with models of spin glasses. The {\it symmetric} version of the present problem 
can be seen to be perfectly equivalent to 
the spherical version of the Sherrington-Kirkpatrick spin-glass with an extra ferromagnetic interaction, on its Nishimori line.
This model was introduced and solved long ago
by a ``spectral method'' using Wigner's semicircle law \cite{PhysRevLett.36.1217}.
Although this analysis is not completely rigorous, it can be made so (hence providing a proof of the replica formula by avoiding the replica trick entirely). 
For the \textit{non-symmetric} inference problem considered in this paper, it is presumably also possible to use a spectral method (using Ginibre's circle law \cite{GinibreLaw}), instead of an interpolation, to arrive at the expression of the mutual information.
However, it has to be noted that the interpolation method presented here
readily extends to rank-one tensor problems.
Indeed, the present analysis can be combined with \cite{BarbierMacrisAllerton2017} to treat the spherical tensors.

We already pointed out that the mutual informations for spherically distributed and i.i.d.\ Gaussian signal vectors are the same.
This is perhaps not so surprising since, roughly speaking, a standard Gaussian vector concentrates on a sphere.
However, we believe that this quick argument is rather weak for two reasons.
First, the only method we know to check that the mutual informations are equal is to compute them separately and compare them.
Secondly, this argument fails when naively applied to the spherical spin-glass model of statistical mechanics.
It is well-known that the spherical and Gaussian spin-glass models are not equivalent (this is goes back to \cite{PhysRev.86.821}, see \cite{Barra2014,Genovese2015} for interesting recent developments).
From this perspective, it is not obvious that in inference the two distributions lead to identical asymptotic mutual informations.

\section*{Acknowledgment}
C. L. acknowledges funding from the  Swiss National Science Foundation, under grant no 200021E-175541. 


\bibliographystyle{IEEEtran}
\bibliography{IEEEabrv,bibliography}

\appendices
\newpage
\section{Proof of Lemma~\ref{lemma:free_entropy_spherical_vector}}\label{app:free_entropy_spherical_vector}
Let $\bX \sim P_x$ a $n$-dimensional random vector uniformly distributed on the sphere of radius $\sqrt{n}$.
We are interested in the average mutual information between $\bX$ and $\widetilde{\bY} = \sqrt{m}\,\bX + \widetilde{\bZ}$  in the high-dimensional regime, where $\widetilde{Z}_i \iid \cN(0,1)$ for $i=1,\dots, n$ and $m > 0$.
We first link this average mutual information to the free entropy $\widetilde{f}_n \triangleq \frac{1}{n}\E \ln \int dP_x(\bx) e^{-\cH(\bx, \bY)}$ where
$\cH(\bx, \bY) \triangleq \sum_{i=1}^{n}
\frac{m}{2} x_i^2 - \sqrt{m}\, x_i \widetilde{Y}_i$. We have:
\begin{align}
\frac{I(\bX ; \widetilde{\bY})}{n} 
= -\frac{h(\bY)}{n} - \frac{h(\bY \vert \bX)}{n}
= -\frac{\E \ln \int dP_x(\bx) e^{-\cH(\bx, \bY) - \frac{\Vert \widetilde{\bY} \Vert^2}{2}}}{n} + \frac{\E \ln e^{-\frac{\Vert \widetilde{\bZ} \Vert^2}{2}}}{n}
= \frac{m}{2} - \widetilde{f}_n \;.
\end{align}
Therefore, proving Lemma~\ref{lemma:free_entropy_spherical_vector} is equivalent to proving
\begin{equation}\label{free_entropy_noisy_spherical_vector}
\lim_{n \to +\infty} \widetilde{f}_n
= \frac{m}{2} - \frac{\ln(1+m)}{2} \;.
\end{equation}
We use a classical interpolation scheme to prove \eqref{free_entropy_noisy_spherical_vector}.
For $t \in [0,1]$, consider the estimation of the $n$-dimensional vector $\widetilde{\bX}$ whose entries are i.i.d.\ with respect to $\cN(0,1)$ from the observations
\begin{align}\label{interpolation_sphere}
\begin{cases}
\bY^{(t)} &= \sqrt{m (1-t)n}\,\frac{\widetilde{\bX}}{\Vert \widetilde{\bX} \Vert} + \bZ\\
\widetilde{\bY}^{(t)} &= \sqrt{m t} \, \widetilde{\bX} \; + \; \widetilde{\bZ}\\
\end{cases}
\end{align}
where the noises $\bZ \in \R^{n}$, $\widetilde{\bZ} \in \R^{n}$ have i.i.d.\ entries with respect to $\cN(0,1)$.
The associated interpolating Hamiltonian is:
\begin{equation}\label{interpolating_hamiltonian_sphere}
\cH_{t}(\widetilde{\bx}; \bY^{(t)}, \widetilde{\bY}^{(t)})
\triangleq
\sum_{i=1}^{n}
\frac{m(1-t)n}{2} \frac{\widetilde{x}_i^2}{\Vert \widetilde{\bx} \Vert^2} - \sqrt{m(1-t)n}\, \frac{\widetilde{x}_i}{\Vert \widetilde{\bx} \Vert} Y_{i}^{(t)}
+
\sum_{i=1}^{n}
\frac{mt}{2} \widetilde{x}_i^2 - \sqrt{mt}\, \widetilde{x}_i \widetilde{Y}_{i}^{(t)}\;.
\end{equation}
Define the interpolating free entropy $\widetilde{f}_n(t) \triangleq \frac1n \E \ln \cZ_{t}(\bY^{(t)},\widetilde{\bY}^{(t)})$ where 
\begin{equation*}
\cZ_{t}(\bY^{(t)},\widetilde{\bY}^{(t)}) = \int  \frac{d\widetilde{\bx}}{\sqrt{2 \pi}^n}e^{-\frac{\Vert \widetilde{\bx} \Vert^2}{2}} \, e^{-\cH_{t}(\widetilde{\bx} ; \bY^{(t)},\widetilde{\bY}^{(t)})} \;.
\end{equation*}
Note that $\sqrt{n} \frac{\widetilde{\bX}}{\Vert \widetilde{\bX} \Vert}$ has the same distribution than $\bX$, i.e., it is uniformly distributed on the $(n-1)$-sphere of radius $\sqrt{n}$.
Then the observation model at $t=0$ is identical to \eqref{spherical_vector_estimation} and we find that $\widetilde{f}_n(0) = \widetilde{f}_n$ is the free entropy whose limit we want to compute.
At $t=1$, the above integral is a simple Gaussian integral and we find $\widetilde{f}_n(1) = \frac{m}{2} - \frac{\ln(1+m)}{2}$. Hence, we have:
\begin{equation}\label{upperbound_difference_f_n}
 \bigg\vert \frac{m}{2} - \frac{\ln(1+m)}{2} - \widetilde{f}_n \bigg\vert
= \bigg\vert \int_0^1 \widetilde{f}'_n(t) \,dt \bigg\vert
\leq \int_0^1 \big\vert \widetilde{f}'_n(t)\big\vert \,dt \;.
\end{equation}
Computing $\widetilde{f}'_n(t)$ is done much like in the proof of Lemma~\ref{lemma:formula_derivative_free_entropy} when computing the derivative of the average free entropy \eqref{interpolating_free_entropy}.
We obtain:
\begin{equation}
\widetilde{f}'_n(t) = \frac{m}{2} \,\E\,\bigg\langle \frac{\widetilde{\bx}^{\sT}\widetilde{\bX}}{\Vert \widetilde{\bx} \Vert \Vert \widetilde{\bX} \Vert}\bigg(\frac{\Vert \widetilde{\bx} \Vert \Vert \widetilde{\bX} \Vert}{n} - 1\bigg)\!\bigg\rangle_{\! t} \;,
\end{equation}
where the angular brackets $\langle - \rangle_t$ denote the expectation w.r.t.\ the posterior distribution
$$
dP(\widetilde{\bx} \vert \bY^{(t)},\widetilde{\bY}^{(t)} )
= \frac{1}{\cZ_{t}(\bY^{(t)},\widetilde{\bY}^{(t)})}\frac{d\widetilde{\bx}}{\sqrt{2 \pi}^n}e^{-\frac{\Vert \widetilde{\bx} \Vert^2}{2}} \, e^{-\cH_{t}(\widetilde{\bx} ; \bY^{(t)},\widetilde{\bY}^{(t)})} \;.
$$
We split $\widetilde{f}'_n(t)$ in two pieces:
$$
f'_n(t)
=
\frac{m}{2} \E\,\bigg\langle \frac{\widetilde{\bx}^{\sT}\widetilde{\bX}}{\Vert \widetilde{\bx} \Vert \Vert \widetilde{\bX} \Vert}
\frac{\Vert \widetilde{\bx} \Vert}{\sqrt{n}}\bigg(\frac{\Vert \widetilde{\bX} \Vert}{\sqrt{n}} - 1\bigg)\!\bigg\rangle_{\! t}
+
\frac{m}{2} \E\,\bigg\langle \frac{\widetilde{\bx}^{\sT}\widetilde{\bX}}{\Vert \widetilde{\bx} \Vert \Vert \widetilde{\bX} \Vert}\bigg(\frac{\Vert \widetilde{\bx} \Vert}{\sqrt{n}} - 1\bigg)\!\bigg\rangle_{\! t} \,.
$$
Applying Cauchy-Schwarz inequality separately to these two pieces, we get:
\begin{align}
\big\vert \widetilde{f}'_n(t) \big\vert
&\leq
\frac{m}{2} \sqrt{\E\,\bigg\langle \frac{(\widetilde{\bx}^{\sT}\widetilde{\bX})^2}{\Vert \widetilde{\bx} \Vert^2 \Vert \widetilde{\bX} \Vert^2}
\frac{\Vert \widetilde{\bx} \Vert^2}{n}\,\bigg\rangle_{\! t}
\E\bigg[\bigg(\frac{\Vert \widetilde{\bX} \Vert}{\sqrt{n}} - 1\bigg)^{\!\! 2}\,\bigg]}
+
\frac{m}{2}
\sqrt{\E\,\bigg\langle \frac{(\widetilde{\bx}^{\sT}\widetilde{\bX})^2}{\Vert \widetilde{\bx} \Vert^2 \Vert \widetilde{\bX} \Vert^2}\!\bigg\rangle_{\! t}
\E\,\bigg\langle \bigg(\frac{\Vert \widetilde{\bx} \Vert}{\sqrt{n}} - 1\bigg)^{\!\! 2} \,\bigg\rangle_{\! t}}\nonumber\\
&\leq
\frac{m}{2} \sqrt{\E\,\bigg\langle \frac{\Vert \widetilde{\bx} \Vert^2}{n}\,\bigg\rangle_{\! t}
	\E\bigg[\bigg(\frac{\Vert \widetilde{\bX} \Vert}{\sqrt{n}} - 1\bigg)^{\!\! 2}\,\bigg]}
+ \frac{m}{2}
\sqrt{\E\,\bigg\langle \bigg(\frac{\Vert \widetilde{\bx} \Vert}{\sqrt{n}} - 1\bigg)^{\!\! 2} \,\bigg\rangle_{\! t}}\nonumber\\
&= m \sqrt{\E\bigg[\bigg(\frac{\Vert \widetilde{\bX} \Vert}{\sqrt{n}} - 1\bigg)^{\!\! 2}\,\bigg]} \;.\label{upperbound_derivative_f_n}
\end{align}
The second inequality follows again from the Cauchy-Schwarz inequality: $\vert \widetilde{\bx}^{\sT}\widetilde{\bX} \vert \leq \Vert \widetilde{\bx} \Vert \Vert \widetilde{\bX} \Vert$.
The subsequent equality is an application of the Nishimori identity (see Lemma~\ref{lemma:nishimori} directly following the proof):
$\E\,\big\langle \frac{\Vert \widetilde{\bx} \Vert^2}{n}\big\rangle_{t} = \E\,\frac{\Vert \widetilde{\bX} \Vert^2}{n}$ and 
$\E\,\big\langle \big(\frac{\Vert \widetilde{\bx} \Vert}{\sqrt{n}} - 1\big)^{2}\,\big\rangle_{t}
=\E\,\big(\frac{\Vert \widetilde{\bX} \Vert}{\sqrt{n}} - 1\big)^{2}$.
The upper bound \eqref{upperbound_derivative_f_n} on the absolute value of the derivative of the interpolating free entropy is valid for all $t \in [0,1]$.
Plugging it back in \eqref{upperbound_difference_f_n} yields:
\begin{equation}\label{final_upperbound_difference_free_entropy_limit}
\bigg\vert \frac{m}{2} - \frac{\ln(1+m)}{2} - \widetilde{f}_n \bigg\vert
\leq \frac{m}{\sqrt{n}} \sqrt{\E\big[\big(\Vert \widetilde{\bX} \Vert - \sqrt{n}\,\big)^{2}\,\big]} \;.
\end{equation}
There exists a constant $C$ such that $\mathbb{P}(\vert \Vert \widetilde{\bX} \Vert - \sqrt{n}\vert \geq a) \leq 2e^{-Ca^2}$ for all $a \geq 0$ (see \cite[Theorem 3.1.1]{vershynin_2018}).
This directly implies $\E\big[\big(\Vert \widetilde{\bX} \Vert - \sqrt{n}\,\big)^{2}\big] \leq \nicefrac{2}{C}$.
Given the upper bound \eqref{final_upperbound_difference_free_entropy_limit}, it concludes the proof of \eqref{free_entropy_noisy_spherical_vector}. \hfill\IEEEQED

\begin{lemma}[Nishimori identity] \label{lemma:nishimori}
Let $(\bX,\bY) \in \R^{n_1} \times \R^{n_2}$ be a pair of jointly distributed random vectors.
Let $k \geq 1$.
Let $\bX^{(1)}, \dots, \bX^{(k)}$ be $k$ independent samples drawn from the conditional distribution $P(\bX=\cdot\, \vert \bY)$, independently of every other random variables.
The angular brackets $\langle - \rangle$ denote the expectation operator with respect to $P(\bX= \cdot\, \vert \bY)$, while $\E$ denotes the expectation with respect to $(\bX,\bY)$.
Then, for all continuous bounded function $g$ we have:
\begin{equation*}
\E\,\langle g(\bY,\bX^{(1)}, \dots, \bX^{(k)}) \rangle
= \E\,\langle g(\bY,\bX^{(1)}, \dots, \bX^{(k-1)}, \bX) \rangle\;.
\end{equation*}
\end{lemma}
\begin{IEEEproof}
This is a simple consequence of Bayes' formula.
It is equivalent to sample the pair $(\bX,\bY)$ according to its joint distribution, or to first sample $\bY$ according to its marginal distribution and to then sample $\bX$ conditionally to $\bY$ from its conditional distribution $P(\bX=\cdot\,|\bY)$.
Hence the $(k+1)$-tuple $(\bY,\bX^{(1)}, \dots,\bX^{(k)})$ is equal in law to $(\bY,\bX^{(1)},\dots,\bX^{(k-1)},\bX)$.
\end{IEEEproof}
\newpage
\section{Establishing the sum-rule of Proposition~\ref{prop:sum_rule}}\label{app:sum_rule}
Remember that we fixed $\lambda=1$ and this without loss of generality.
We remind the reader of the definitions of the scalar overlaps:
$Q_u \triangleq \frac{1}{n_u} \sum_i u_i U_i$ and
$Q_v \triangleq \frac{1}{n_v} \sum_i v_i V_i$.
\begin{lemma}[Average interpolating free entropy at $t=0$ and $t=1$]
Assume that both $R_u(t,\epsilon)$ and $R_v(t,\epsilon)$ are uniformly bounded in $(t,\epsilon) \in [0,1] \times [0,+\infty)^2$.
The average interpolating free entropy $f_n(t,\epsilon)$ whose definition is given by \eqref{interpolating_free_entropy} satisfies:
\begin{align}
f_n(0,\epsilon) &= f_n(0,0) + \mathcal{O}(\Vert \epsilon \Vert) \;;\label{f_n(0,epsilon)}\\
f_n(1,\epsilon) &=  \alpha_u \varphi(\alpha_v \rho_u R_v(1,\epsilon)) + \alpha_v \varphi(\alpha_u \rho_v R_u(1,\epsilon)) + \smallO_n(1) \;;\label{f_n(1,epsilon)}
\end{align}
where $\smallO_n(1)$ is a quantity that vanishes uniformly in $\epsilon$ as $n$ gets large,
and $\mathcal{O}(\Vert \epsilon \Vert)$ is a quantity whose absolute value is upper bounded by $C \Vert \epsilon \Vert$ for some constant $C$ independent of both $n$ and $\epsilon$.
\end{lemma}
\begin{IEEEproof}
By definition, $f_n(0,\epsilon) = \frac{1}{n} \E\ln \cZ_{0,\epsilon}(\bY^{(0)},\widetilde{\bY}^{(0,\epsilon)}, \overline{\bY}^{(0,\epsilon)})$ where
$$
\cZ_{0,\epsilon}(\bY^{(0)},\widetilde{\bY}^{(0,\epsilon)}, \overline{\bY}^{(0,\epsilon)}) \triangleq \int dP_{u}(\bu)dP_v(\bv) \, e^{-\cH_{0,\epsilon}(\bu, \bv ; \bY^{(0)},\widetilde{\bY}^{(0,\epsilon)}, \overline{\bY}^{(0,\epsilon)})}
$$
and $\cH_{0,\epsilon}(\bu, \bv ; \bY^{(0)},\widetilde{\bY}^{(0,\epsilon)}, \overline{\bY}^{(0,\epsilon)})$ is the Hamiltonian \eqref{interpolating_hamiltonian} evaluated at $t=0$. Remembering that $R_u(0,\epsilon) = \epsilon_u, R_v(0,\epsilon) = \epsilon_v$, and replacing $\bY^{(0)}, \widetilde{\bY}^{(0,\epsilon)}, \overline{\bY}^{(0,\epsilon)}$ by their expressions on the right-hand side of \eqref{interpolation_model}, we obtain:
\begin{equation}\label{expression_Z_0,epsilon}
\cZ_{0,\epsilon}(\bY^{(0)},\widetilde{\bY}^{(0,\epsilon)}, \overline{\bY}^{(0,\epsilon)})
= \int dP_{u}(\bu)dP_v(\bv) \, e^{-\cH_{0,\epsilon}(\bu, \bv ; \bU, \bV, \bZ,\widetilde{\bZ}, \overline{\bZ})}
\end{equation}
with
\begin{align}
\cH_{0,\epsilon}(\bu, \bv ; \bU, \bV, \bZ,\widetilde{\bZ}, \overline{\bZ})
\triangleq
\sum_{i=1}^{n_u}\sum_{j=1}^{n_v}
\frac{u_i^2 v_j^2}{2n}  - \frac{1}{n}\, u_i U_i v_j V_j - \frac{u_i v_j Z_{ij}}{\sqrt{n}}
&+
\sum_{i=1}^{n_u} \frac{\alpha_v \epsilon_v}{2} u_i^2 - \alpha_v \epsilon_v\, u_i U_i - \sqrt{\alpha_v \epsilon_v}\, u_i \widetilde{Z}_i \nonumber\\
&+
\sum_{j=1}^{n_v} \frac{\alpha_u \epsilon_u}{2} v_j^2 - \alpha_u \epsilon_u\, v_j V_j - \sqrt{\alpha_u \epsilon_u}\, v_j \overline{Z}_j \;.
\end{align}
Making use of \eqref{expression_Z_0,epsilon}, the partial derivative of $\epsilon \mapsto f_n(0,\epsilon)$ with respect to $\epsilon_u$ reads:
\begin{align}
\frac{\partial f_n}{\partial \epsilon_u}\bigg\vert_{t=0,\epsilon}
&= -\frac{1}{n}\E\bigg[\frac{\int dP_{u}(\bu)dP_v(\bv) \, \frac{\partial \cH_{0,\epsilon}(\bu, \bv ; \bU, \bV, \bZ,\widetilde{\bZ}, \overline{\bZ})}{\partial \epsilon_u}e^{-\cH_{0,\epsilon}(\bu, \bv ; \bU, \bV, \bZ,\widetilde{\bZ}, \overline{\bZ})}}{\int dP_{u}(\bu)dP_v(\bv) \, e^{-\cH_{0,\epsilon}(\bu, \bv ; \bU, \bV, \bZ,\widetilde{\bZ}, \overline{\bZ})}}\bigg]\nonumber\\
&= -\frac{1}{n}\E\,\bigg\langle \frac{\partial \cH_{0,\epsilon}(\bu, \bv ; \bU, \bV, \bZ,\widetilde{\bZ}, \overline{\bZ})}{\partial \epsilon_u}\bigg\rangle_{\!\! t=0,\epsilon}\nonumber\\
&= - \frac{\alpha_u}{2n}\sum_{j=1}^{n_v} \E\,\langle   v_j^2 \rangle_{0,\epsilon}
+ \frac{\alpha_u}{n}\sum_{j=1}^{n_v} \E\,\langle   v_j V_j \rangle_{0,\epsilon}
+ \frac{1}{2n}\sqrt{\frac{\alpha_u}{\epsilon_u}}\sum_{j=1}^{n_v}\E\,\langle v_j \overline{Z}_j \rangle_{0,\epsilon} \;. \label{d_fn/d_epsilon_u}
\end{align}
We can now simplify the expectation $\E\,\langle v_j \overline{Z}_j \rangle_{0,\epsilon}$ with an integration by parts with respect to the standard Gaussian $\overline{Z}_j$:
\begin{align}
\E\,\langle v_j \overline{Z}_j \rangle_{0,\epsilon}
&= \E\bigg[\overline{Z}_j \frac{\int dP_{u}(\bu)dP_v(\bv) \, v_je^{-\cH_{0,\epsilon}(\bu, \bv ; \bU, \bV, \bZ,\widetilde{\bZ}, \overline{\bZ})}}{\int dP_{u}(\bu)dP_v(\bv) \, e^{-\cH_{0,\epsilon}(\bu, \bv ; \bU, \bV, \bZ,\widetilde{\bZ}, \overline{\bZ})}}\bigg]\nonumber\\
&= \E\bigg[\frac{\int dP_{u}(\bu)dP_v(\bv) \, v_j e^{-\cH_{0,\epsilon}(\bu, \bv ; \bU, \bV, \bZ,\widetilde{\bZ}, \overline{\bZ})}
\int dP_{u}(\bu)dP_v(\bv) \, \frac{\partial \cH_{0,\epsilon}(\bu, \bv ; \bU, \bV, \bZ,\widetilde{\bZ}, \overline{\bZ})}{\partial \overline{Z}_j}e^{-\cH_{0,\epsilon}(\bu, \bv ; \bU, \bV, \bZ,\widetilde{\bZ}, \overline{\bZ})}}{\big(\int dP_{u}(\bu)dP_v(\bv) \, e^{-\cH_{0,\epsilon}(\bu, \bv ; \bU, \bV, \bZ,\widetilde{\bZ}, \overline{\bZ})}\big)^2}\bigg]\nonumber\\
&\qquad-\E\bigg[\frac{\int dP_{u}(\bu)dP_v(\bv) \, v_j \frac{\partial \cH_{0,\epsilon}(\bu, \bv ; \bU, \bV, \bZ,\widetilde{\bZ}, \overline{\bZ})}{\partial \overline{Z}_j}e^{-\cH_{0,\epsilon}(\bu, \bv ; \bU, \bV, \bZ,\widetilde{\bZ}, \overline{\bZ})}}{\int dP_{u}(\bu)dP_v(\bv) \, e^{-\cH_{0,\epsilon}(\bu, \bv ; \bU, \bV, \bZ,\widetilde{\bZ}, \overline{\bZ})}}\bigg]\nonumber\\
&= \E\bigg[\langle v_j \rangle_{0,\epsilon} \bigg\langle \frac{\partial \cH_{0,\epsilon}(\bu, \bv ; \bU, \bV, \bZ,\widetilde{\bZ}, \overline{\bZ})}{\partial \overline{Z}_j} \bigg\rangle_{\!\! 0,\epsilon}\,\bigg]
-\E\bigg[\bigg\langle v_j \frac{\partial \cH_{0,\epsilon}(\bu, \bv ; \bU, \bV, \bZ,\widetilde{\bZ}, \overline{\bZ})}{\partial \overline{Z}_j}\bigg\rangle_{\!\! 0,\epsilon}\,\bigg]\nonumber\\
&= -\sqrt{\alpha_u \epsilon_u} \,\E[\langle v_j \rangle_{0,\epsilon}^2]
+ \sqrt{\alpha_u \epsilon_u} \,\E[\langle v_j^2\rangle_{0,\epsilon}]\;.\label{E_v_barZ_IBP}
\end{align}
Plugging \eqref{E_v_barZ_IBP} back in \eqref{d_fn/d_epsilon_u} yields:
\begin{equation}
\frac{\partial f_n}{\partial \epsilon_u}\bigg\vert_{t=0,\epsilon}
= \frac{\alpha_u}{n}\sum_{j=1}^{n_v} \E\,\langle   v_j V_j \rangle_{0,\epsilon}
- \frac{\alpha_u}{2n}\sum_{j=1}^{n_v} \E[\langle v_j \rangle_{0,\epsilon}^2]
= \frac{\alpha_u}{2}\frac{n_v}{n} \E\,\langle   Q_v \rangle_{0,\epsilon} \;.
\end{equation}
The second equality follows from the Nishimori identity $\E[\langle v_j \rangle_{0,\epsilon}^2] = \E[\langle v_j \rangle_{0,\epsilon} V_j]$ (see Lemma~\ref{lemma:nishimori}).
We have proved that $\forall \epsilon \in [0,+\infty)^2: \frac{\partial f_n}{\partial \epsilon_u}\big\vert_{t=0,\epsilon} = \frac{\alpha_u}{2}\frac{n_v}{n} \E\,\langle   Q_v \rangle_{0,\epsilon}$ and (this is proved in a similar way) $\frac{\partial f_n}{\partial \epsilon_v}\big\vert_{t=0,\epsilon} = \frac{\alpha_v}{2}\frac{n_u}{n} \E\,\langle   Q_u \rangle_{0,\epsilon}$.
Besides, by Cauchy-Schwarz inequality, $\vert Q_u \vert \leq \nicefrac{\Vert \bu \Vert \Vert \bU \Vert}{n_u} = \rho_u$ and $\vert Q_v \vert \leq \nicefrac{\Vert \bv \Vert \Vert \bV \Vert}{n_v} = \rho_v$ almost surely.
By the mean-value theorem -- $\epsilon = (\epsilon_u,\epsilon_v)$ -- :
\begin{align*}
\vert f_n(0,\epsilon) - f_n(0,0) \vert
&\leq \vert f_n(0,\epsilon) - f_n(0,(0,\epsilon_v)) \vert  + \vert f_n(0,(0,\epsilon_v)) - f_n(0,0) \vert\\
&\leq \bigg\vert \frac{\partial f_n}{\partial \epsilon_u}\Big\vert_{0,\epsilon} \bigg\vert \vert \epsilon_u \vert
+ \bigg\vert \frac{\partial f_n}{\partial \epsilon_v}\Big\vert_{0,(0,\epsilon_v)} \bigg\vert \vert \epsilon_v \vert\\
&\leq \frac{\alpha_u \rho_u}{2}\frac{n_v}{n}\vert \epsilon_u \vert
+ \frac{\alpha_v \rho_v}{2}\frac{n_u}{n} \vert \epsilon_v \vert \;.
\end{align*}
Knowing that $(\nicefrac{n_u}{n}, \nicefrac{n_v}{n}) \to (\alpha_u,\alpha_v)$, this last upper bound concludes the proof of \eqref{f_n(0,epsilon)}.

At $t=1$, the observation $\bY^{(1)} = \bZ$ is pure noise while $\widetilde{\bY}^{(1,\epsilon)} = \sqrt{\alpha_v R_v(1,\epsilon)}\, \bU + \widetilde{\bZ}$ and $\overline{\bY}^{(1,\epsilon)} = \sqrt{\alpha_u R_u(1,\epsilon)}\, \bV + \overline{\bZ}$ are two decoupled channels like the one described in Lemma~\ref{lemma:free_entropy_spherical_vector}. Then, we easily see that $f_n(1,\epsilon) = \frac{n_u}{n} \widetilde{f}_{n_u} + \frac{n_v}{n} \widetilde{f}_{n_v}$ where
\begin{align*}
\widetilde{f}_{n_u} &\triangleq \frac{1}{n_u}\E \ln \int dP_u(\bu) \,e^{-\sum_{i=1}^{n_u}
	\frac{\alpha_v R_v(1,\epsilon)}{2} u_i^2 - \sqrt{\alpha_v R_v(1,\epsilon)}\, u_i \widetilde{Y}_i^{(1,\epsilon)}} \;;\\
\widetilde{f}_{n_v} &\triangleq \frac{1}{n_v}\E \ln \int dP_v(\bv) \,e^{-\sum_{i=1}^{n_v}
	\frac{\alpha_u R_u(1,\epsilon)}{2} v_i^2 - \sqrt{\alpha_u R_u(1,\epsilon)}\, v_i \overline{Y}_i^{(1,\epsilon)}} \;;
\end{align*}
are the average free entropy of the two aforementioned channels.
In the proof of Lemma~\ref{lemma:free_entropy_spherical_vector}, we show that
\begin{equation}\label{upperbound_diff_ftilde_limit}
\vert \widetilde{f}_{n_u} - \varphi(\alpha_v \rho_u R_v(1,\epsilon)) \vert \leq \frac{2 \alpha_v R_v(1,\epsilon)}{C \sqrt{n_u}}
\quad\text{and}\quad
\vert \widetilde{f}_{n_v} - \varphi(\alpha_u \rho_v R_u(1,\epsilon)) \vert \leq \frac{2 \alpha_u R_u(1,\epsilon)}{C \sqrt{n_v}}
\end{equation}
where $C$ is a constant independent of $\epsilon$.
The two upper bounds \eqref{upperbound_diff_ftilde_limit} together with the assumption on the uniform boundedness of $R_u,R_v$ yields \eqref{f_n(1,epsilon)}:
$$
f_n(1,\epsilon) = \frac{n_u}{n} \widetilde{f}_{n_u} + \frac{n_v}{n} \widetilde{f}_{n_v}
= \alpha_u \varphi(\alpha_v \rho_u R_v(1,\epsilon)) + \alpha_v \varphi(\alpha_u \rho_v R_u(1,\epsilon)) + \smallO_n(1) \;.
$$
\end{IEEEproof}
\begin{lemma}[Derivative of the average interpolating free entropy]\label{lemma:formula_derivative_free_entropy}
Denote $R'_u(\cdot,\epsilon)$ and $R'_v(\cdot,\epsilon)$ the derivative of $R_u(\cdot,\epsilon)$ and $R_v(\cdot,\epsilon)$, respectively.
Assume that both $R'_u(t,\epsilon)$ and $R'_v(t,\epsilon)$ are uniformly bounded in $(t,\epsilon)$ belonging to $[0,1] \times [0,+\infty)^2$.
Then, the derivative with respect to $t$ of the average free entropy \eqref{interpolating_free_entropy} satisfies $\forall (t, \epsilon) \in [0,1] \times [0,+\infty)^2$:
\begin{equation}\label{formula_derivative_free_entropy}
f'_{n}(t,\epsilon)
= -\frac{\alpha_u \alpha_v}{2}  
\E\big[\big\langle \big( Q_u - R'_u(t,\epsilon)\big)\big(Q_v - R'_v(t,\epsilon)\big) \big\rangle_{t,\epsilon} \,\big]
+\frac{\alpha_u \alpha_v}{2} R'_u(t,\epsilon) R'_v(t,\epsilon)
+ \smallO_{n}(1)\;,
\end{equation}
where $\smallO_{n}(1)$ vanishes uniformly in $(t, \epsilon)$ as $n$ goes to infinity.
\end{lemma}
\begin{IEEEproof}
The conditional probability density function of $(\bY^{(t)},\widetilde{\bY}^{(t,\epsilon)} , \overline{\bY}^{(t,\epsilon)})$ given $(\bU,\bV)$ is
$$
p_{\bY^{(t)},\widetilde{\bY}^{(t,\epsilon)} , \overline{\bY}^{(t,\epsilon)} \vert \bU,\bV}(\by, \widetilde{\by}, \bar{\by} \vert \bu, \bv)
\triangleq \frac{1}{\sqrt{2\pi}^{n_u n_v + n_u + n_v}} e^{-\frac{1}{2} \big(\sum_{i,j} y_{ij}^2+ \Vert \widetilde{\by} \Vert^2 +\Vert \bar{\by} \Vert^2 \big)-\cH_{t,\epsilon}(\bu, \bv; \by,\widetilde{\by}, \bar{\by})}
$$
where
\begin{multline}\label{interpolating_hamiltonian_small_y}
\cH_{t,\epsilon}(\bu, \bv ; \by, \widetilde{\by}, \bar{\by} )
\triangleq
\sum_{i=1}^{n_u}\sum_{j=1}^{n_v}
\frac{1-t}{2n} u_i^2 v_j^2 - \sqrt{\frac{1-t}{n}}\, u_i v_j y_{ij}
+
\sum_{i=1}^{n_u} \frac{\alpha_v R_v(t,\epsilon)}{2} u_i^2 - \sqrt{\alpha_v R_v(t,\epsilon)}\, u_i \widetilde{y}_i\\
+
\sum_{j=1}^{n_v} \frac{\alpha_u R_u(t,\epsilon)}{2} v_j^2 - \sqrt{\alpha_u R_u(t,\epsilon)}\, v_j \overline{y}_j \;.
\end{multline}
Therefore, the average interpolating free entropy \eqref{interpolating_free_entropy} satisfies:
\begin{align}
f_n(t,\epsilon)
&=  \frac{1}{n} \E_{\bU,\bV}\bigg[\E\Big[
\ln \cZ_{t,\epsilon}\big( \bY^{(t)},\widetilde{\bY}^{(t,\epsilon)} , \overline{\bY}^{(t,\epsilon)}\big) \Big\vert \bU, \bV \Big]\bigg] \nonumber\\
&= \frac{1}{n} \E_{\bU,\bV}\bigg[\int d\by d\widetilde{\by} d\bar{\by}\:
\frac{e^{-\frac{1}{2} \big(\sum_{i,j} y_{ij}^2+ \Vert \widetilde{\by} \Vert^2 +\Vert \bar{\by} \Vert^2 \big)}}{\sqrt{2\pi}^{n_u n_v + n_u + n_v}} e^{-\cH_{t,\epsilon}(\bU, \bV; \by,\widetilde{\by}, \bar{\by})}
\ln \cZ_{t,\epsilon}\big(\by,\widetilde{\by}, \bar{\by}\big) \bigg]\;.\label{eq:precise_formula_f}
\end{align}
where we remind that $\cZ_{t,\epsilon}(\by,\widetilde{\by}, \bar{\by}) \triangleq \int dP_{u}(\bu)dP_v(\bv) \, e^{-\cH_{t,\epsilon}(\bu, \bv ; \by,\widetilde{\by}, \bar{\by})}$
Taking the derivative of \eqref{eq:precise_formula_f} with respect to $t$, we directly obtain:
\begin{equation*}
f'_n(t,\epsilon)
= -\frac{1}{n} \underbrace{\E\big[\cH'_{t,\epsilon}\big(\bU, \bV; \bY^{(t)},\widetilde{\bY}^{(t,\epsilon)} , \overline{\bY}^{(t,\epsilon)}\big)\ln \cZ_{t,\epsilon}\big(\bY^{(t)},\widetilde{\bY}^{(t,\epsilon)} , \overline{\bY}^{(t,\epsilon)}\big)\big]}_{\triangleq T_1}
-\frac{1}{n} \underbrace{\E\,\big\langle \cH'_{t,\epsilon}\big(\bu, \bv; \bY^{(t)},\widetilde{\bY}^{(t,\epsilon)} , \overline{\bY}^{(t,\epsilon)}\big)\big\rangle_{\!t,\epsilon}}_{\triangleq T_2}\;,
\end{equation*}
with
\begin{multline}\label{derivative_hamiltonian}
\cH'_{t,\epsilon}(\bu, \bv ; \by,\widetilde{\by}, \overline{\by})
\triangleq
\sum_{i=1}^{n_u}\sum_{j=1}^{n_v}
-\frac{u_i^2 v_j^2}{2n}  + \frac{u_i v_j y_{ij}^{(t)}}{2\sqrt{n(1-t)}} 
+
\sum_{i=1}^{n_u} \frac{\alpha_v R'_v(t,\epsilon)}{2} u_i^2 - \frac{R'_v(t,\epsilon)}{2}\sqrt{\frac{\alpha_v}{R_v(t,\epsilon)}}\, u_i \widetilde{y}_i^{(t,\epsilon)} \\
+
\sum_{j=1}^{n_v} \frac{\alpha_u R'_u(t,\epsilon)}{2} v_j^2 - \frac{R'_u(t,\epsilon)}{2}\sqrt{\frac{\alpha_u}{R_u(t,\epsilon)}}\, v_j \overline{y}_j^{(t,\epsilon)} \;.
\end{multline}
Equation \eqref{derivative_hamiltonian} is obtained by differentiating the interpolating Hamiltonian \eqref{interpolating_hamiltonian_small_y} with respect to $t$.
If we evaluate \eqref{derivative_hamiltonian} at
$(\bu,\bv, \by, \widetilde{\by}, \overline{\by}) = (\bU,\bV,\bY^{(t)},\widetilde{\bY}^{(t,\epsilon)}, \overline{\bY}^{(t,\epsilon)})$, we get:
\begin{multline}\label{eq:H_timederivative}
\cH'_{t,\epsilon}\big(\bU, \bV; \bY^{(t)},\widetilde{\bY}^{(t,\epsilon)} , \overline{\bY}^{(t,\epsilon)}\big) =
\sum_{i=1}^{n_u}\sum_{j=1}^{n_v} \frac{U_i V_j Z_{ij}}{2\sqrt{n(1-t)}}
-
\sum_{i=1}^{n_u} \frac{R'_v(t,\epsilon)}{2}\sqrt{\frac{\alpha_v}{R_v(t,\epsilon)}}\, U_i \widetilde{Z}_i\\
-
\sum_{j=1}^{n_v} \frac{R'_u(t,\epsilon)}{2}\sqrt{\frac{\alpha_u}{R_u(t,\epsilon)}}\, V_j \overline{Z}_j \;.
\end{multline}
$T_2$ is now easily shown to be zero thanks to the Nishimori identity (see Lemma~\ref{lemma:nishimori}):
\begin{align*}
T_2 &= \E\,\big\langle \cH'_{t,\epsilon}(\bu, \bv; \bY^{(t)},\widetilde{\bY}^{(t,\epsilon)} , \overline{\bY}^{(t,\epsilon)}\big) \big\rangle_{t,\epsilon}\\
&= \E\big[ \cH'_{t,\epsilon}\big(\bU, \bV; \bY^{(t)},\widetilde{\bY}^{(t,\epsilon)} , \overline{\bY}^{(t,\epsilon)}\big)\big]\\
&=\sum_{i=1}^{n_u}\sum_{j=1}^{n_v} \frac{\E[U_i V_j Z_{ij}]}{2\sqrt{n(1-t)}}
-
\sum_{i=1}^{n_u} \frac{R'_v(t,\epsilon)}{2}\sqrt{\frac{\alpha_v}{R_v(t,\epsilon)}}\, \E[U_i \widetilde{Z}_i]
-
\sum_{j=1}^{n_v} \frac{R'_u(t,\epsilon)}{2}\sqrt{\frac{\alpha_u}{R_u(t,\epsilon)}}\, \E[V_j \overline{Z}_j]\\
&=0\:.
\end{align*}
Therefore, $f'_n(t,\epsilon) = \nicefrac{-T_1}{n}$. Plugging~\eqref{eq:H_timederivative} in the definition of $T_1$, we obtain:
\begin{multline}\label{eq:non_final_formula_derivative}
f'_{n}(t,\epsilon)
= -\sum_{i=1}^{n_u}\sum_{j=1}^{n_v} \frac{\E[ U_i V_j Z_{ij} \ln \cZ_{t,\epsilon}(\bY^{(t)},\widetilde{\bY}^{(t,\epsilon)} , \overline{\bY}^{(t,\epsilon)})]}{2 n \sqrt{n(1-t)}}
+
\sum_{i=1}^{n_u} \frac{R'_v(t,\epsilon)}{2 n}\sqrt{\frac{\alpha_v}{R_v(t,\epsilon)}}\, \E[U_i \widetilde{Z}_i\ln \cZ_{t,\epsilon}(\bY^{(t)},\widetilde{\bY}^{(t,\epsilon)} , \overline{\bY}^{(t,\epsilon)})]\\
+
\sum_{j=1}^{n_v} \frac{R'_u(t,\epsilon)}{2 n}\sqrt{\frac{\alpha_u}{R_u(t,\epsilon)}}\, \E[V_j \overline{Z}_j \ln \cZ_{t,\epsilon}(\bY^{(t)},\widetilde{\bY}^{(t,\epsilon)} , \overline{\bY}^{(t,\epsilon)})]\;.
\end{multline}
The three expectations appearing on the right-hand side of~\eqref{eq:non_final_formula_derivative} are simplified thanks to Stein's lemma, i.e., by integrating by parts w.r.t.\ the Gaussian noises:
\begin{align*}
\E[ U_i V_j Z_{ij} \ln \cZ_{t,\epsilon}(\bY^{(t)},\widetilde{\bY}^{(t,\epsilon)} , \overline{\bY}^{(t,\epsilon)})]
&= \E\bigg[ U_i V_j \frac{\partial \ln \cZ_{t,\epsilon}(\bY^{(t)},\widetilde{\bY}^{(t,\epsilon)} , \overline{\bY}^{(t,\epsilon)})}{\partial Z_{ij}}\bigg]\\
&= -\E\bigg[ U_i V_j \bigg\langle \frac{\partial \cH_{t,\epsilon}(\bu,\bv ; \bY^{(t)},\widetilde{\bY}^{(t,\epsilon)} , \overline{\bY}^{(t,\epsilon)})}{\partial Z_{ij}}\bigg\rangle_{\!\! t, \epsilon}\,\bigg]\\
&= \sqrt{\frac{1-t}{n}}\E\,\langle u_i U_i v_j V_j \rangle_{t,\epsilon} \;.
\end{align*}
In a similar way:
\begin{align*}
\E[ U_i \widetilde{Z}_{i} \ln \cZ_{t,\epsilon}(\bY^{(t)},\widetilde{\bY}^{(t,\epsilon)} , \overline{\bY}^{(t,\epsilon)})]
&= \sqrt{\alpha_v R_v(t,\epsilon)}\E\,\langle u_i U_i \rangle_{t,\epsilon} \:;\\
\E[ V_j \overline{Z}_{j} \ln \cZ_{t,\epsilon}(\bY^{(t)},\widetilde{\bY}^{(t,\epsilon)} , \overline{\bY}^{(t,\epsilon)})]
&= \sqrt{\alpha_u R_u(t,\epsilon)}\E\,\langle v_j V_j \rangle_{t,\epsilon} \:.
\end{align*}
Hence, we have:
\begin{equation}
f'_{n}(t,\epsilon)
=
- \frac{1}{2}\frac{n_u}{n}\frac{n_v}{n}\E\,\langle Q_u Q_v \rangle_{t,\epsilon}
+ \frac{n_u}{n}\frac{\alpha_v R'_v(t,\epsilon)}{2}\E\,\langle Q_u \rangle_{t,\epsilon}
+ \frac{n_v}{n}\frac{\alpha_u R'_u(t,\epsilon)}{2}\E\,\langle Q_v \rangle_{t,\epsilon} \;.
\end{equation}
By Cauchy-Schwarz inequality, $\vert Q_u \vert \leq \nicefrac{\Vert \bu \Vert \Vert \bU \Vert}{n_u} = \rho_u$ and $\vert Q_v \vert \leq \nicefrac{\Vert \bv \Vert \Vert \bV \Vert}{n_v} = \rho_v$ almost surely.
It comes:
\begin{align}
f'_{n}(t,\epsilon)
=
- \frac{\alpha_u \alpha_v}{2} \E\,\langle (Q_u - R'_u(t,\epsilon)) (Q_v - R'_v(t,\epsilon)) \rangle_{t,\epsilon}
+ \frac{\alpha_u \alpha_v}{2} R'_u(t,\epsilon) R'_v(t,\epsilon)
+ \smallO_n(1) \:,
\end{align}
where $\smallO_n(1)$ is a quantity that vanishes uniformly in $(t, \epsilon) \in [0,1] \times [0,+\infty)^2$ as $n$ goes to infinity.
\end{IEEEproof}
\newpage
\section{Properties of the solutions to the ordinary differential equations governing the interpolation paths}\label{app:properties_ode}
This appendix is dedicated to the proof of Proposition~\ref{prop:ode} in Section~\ref{section:proof_adaptive_interpolation} of the main text.
We first recall a few definitions.
For $t \in [0,1]$ and $R = (R_u, R_v) \in [0, +\infty)^2$, consider the problem of estimating $(\bU, \bV)$ from the observations:
\begin{align*}
\begin{cases}
\bY^{(t)} &= \sqrt{\frac{1-t}{n}}\,\bU \, \bV^{\sT} + \bZ\;;\\
\widetilde{\bY}^{(t,R_v)} &= \sqrt{\alpha_v R_v}\, \bU + \widetilde{\bZ} \;\;\:;\\
\overline{\bY}^{(t,R_u)} &= \sqrt{\alpha_u R_u}\, \bV + \overline{\bZ}\;\;;
\end{cases}
\end{align*}
where $\bU \sim P_u$, $\bV \sim P_v$ and all of the noises $\bZ \in \R^{n_u \times n_v}$, $\widetilde{\bZ} \in \R^{n_u}$, $\overline{\bZ} \in \R^{n_v}$ have i.i.d.\ entries with respect to $\cN(0,1)$.
The posterior distribution of $(\bU, \bV)$ given $(\bY^{(t)},\widetilde{\bY}^{(t,R_v)}, \overline{\bY}^{(t,R_u)})$ reads
\begin{equation*}
dP(\bu, \bv \vert \bY^{(t)},\widetilde{\bY}^{(t,R_v)}, \overline{\bY}^{(t,R_u)})
= \frac{dP_{u}(\bu)dP_v(\bv) \, e^{-\cH_{t,R}(\bu, \bv ; \bY^{(t)},\widetilde{\bY}^{(t,R_v)}, \overline{\bY}^{(t,R_u)})}}{ \cZ_{t,R}(\bY^{(t)},\widetilde{\bY}^{(t,R_v)}, \overline{\bY}^{(t,R_u)})} \;,
\end{equation*}
where $\cZ_{t,R}(\bY^{(t)},\widetilde{\bY}^{(t,R_v)}, \overline{\bY}^{(t,R_u)})$ is the normalization factor and $\cH_{t,R}$ denotes the associated interpolating Hamiltonian:
\begin{align*}
\cH_{t,R}(\bu, \bv ; \bY^{(t)}, \widetilde{\bY}^{(t,R_v)}, \overline{\bY}^{(t,R_u)})
\triangleq
\sum_{i=1}^{n_u}\sum_{j=1}^{n_v}
\frac{1-t}{2n} u_i^2 v_j^2 - \sqrt{\frac{1-t}{n}}\, u_i v_j Y_{ij}^{(t)}
&+ \sum_{i=1}^{n_u} \frac{\alpha_v R_v}{2} u_i^2 - \sqrt{\alpha_v R_v}\, u_i \widetilde{Y}_i^{(t,R_v)}\\
&+ \sum_{j=1}^{n_v} \frac{\alpha_u R_u}{2} v_j^2 - \sqrt{\alpha_u R_u}\, v_j \overline{Y}_j^{(t,R_u)} \;.
\end{align*}
The angular brackets $\langle - \rangle_{t,R}$ denote the expectation w.r.t.\ this last the posterior.
Let $m_u \in [0,\rho_u]$, $F_v(t, R) \triangleq \E \langle Q_v \rangle_{t,R}$ and $F_u(t, R) \triangleq 2 \rho_u \varphi^{\prime}(\alpha_v \rho_u \E \langle Q_v \rangle_{t,R})$.
We consider the two following second-order ODEs with initial value $\epsilon \in [0,+\infty)^2$:
\begin{align*}
g' = \;\big(m_u , F_v(t,g)\big)\quad\;\;\, , \: g(0)=\epsilon \;;\\
g' = \big(F_u(t,g), F_v(t, g)\big) \, , \: g(0)=\epsilon\;.
\end{align*}
We can now repeat and prove Proposition~\ref{prop:ode}.
\begin{proposition*}
	For all $\epsilon \in [0,+\infty)^2$, there exists a unique global solution $R(\cdot,\epsilon): [0,1] \to [0,+\infty)^2 $ to
	\begin{equation*}
	g' = \big(F_u(t,g), F_v(t, g)\big) \, , \: g(0)=\epsilon\;.
	\end{equation*}
	This solution is continuously differentiable and its derivative $R'(\cdot,\epsilon)$ satisfies $R'([0,1],\epsilon) \subseteq [0, \rho_u] \times [0, \rho_v]$.
	Besides, for all $t \in [0,1]$, $R(t,\cdot)$ is a $\mathcal{C}^1$-diffeomorphism from $[0, +\infty)^2$ into its image whose Jacobian determinant is greater than, or equal to, one:
	\begin{equation*}
	\forall \,\epsilon \in [0, +\infty)^2: \det J_{R(t,\cdot)}(\epsilon) \geq 1 \:,
	\end{equation*}
	where $J_{R(t,\cdot)}$ denotes the Jacobian matrix of $R(t,\cdot)$.\\
	Finally, the same statement holds if instead we consider the second-order ODE:
	\begin{equation*}
	g' = \;\big(m_u , F_v(t,g)\big)\quad\;\;\, , \: g(0)=\epsilon \;.
	\end{equation*}
\end{proposition*}
\begin{IEEEproof}
We only give the proof for the ODE $g' = \big(F_u(t,g), F_v(t, g)\big)$, the one for $g' = \;\big(m_u , F_v(t,g)\big)$ is simpler and follows the same arguments.

By Nishimori identity (see Lemma~\ref{lemma:nishimori}) and Jensen's inequality:
$$
\E\langle Q_v \rangle_{t,R} = \frac{\E\Vert\langle \bv \rangle_{t,R}\Vert^2}{n_v}
\leq \frac{\E\langle \Vert\bv\Vert^2 \rangle_{t,R}}{n_v} = \frac{\E\,\Vert\bV\Vert^2}{n_v} = \rho_v \;,
$$
i.e., $\E\langle Q_v \rangle_{t,R} \in [0,\rho_v]$ for all $(t,R) \in [0,1] \times [0,+\infty)^2$.
Thus, the function $F: (t, R) \mapsto (F_u(t,R), F_v(t,R))$ is defined on all $[0,1] \times [0,+\infty)^2$ and takes value in $[0,\rho_u] \times [0, \rho_v]$.
To invoke Cauchy-Lipschitz theorem, we have to check that $F$ is continuous in $t$ and uniformly Lipschitz continuous in $R$ (meaning the Lipschitz constant is independent of $t$).
The continuity is easy to show using domination arguments.
To check the uniform Lipschitzianity, we show that the Jacobian matrix $J_{F(t,\cdot)}(R)$ of $F(t,\cdot)$ is uniformly bounded in $(t, R)$:
\begin{equation*}
J_{F(t,\cdot)}(R) =
\begin{bmatrix}
c(t,R) & c(t,R)\\
1 & 1
\end{bmatrix}
\begin{bmatrix}
\frac{\partial F_v}{\partial R_u}\big\vert_{t,R} & 0\\
0 & \frac{\partial F_v}{\partial R_v}\big\vert_{t,R}
\end{bmatrix}\:,
\end{equation*}
with $c(t,R) \triangleq 2\alpha_v\rho_u^2\varphi''(\alpha_v \rho_u F_v(t,R)) \in [0,\alpha_v \rho_u^2]$ and
\begin{align*}
\frac{\partial F_v}{\partial R_u}
&= \frac{\alpha_v}{n_v} \sum_{i=1}^{n_v} \sum_{j=1}^{n_v} \E[(\langle v_i v_j \rangle_{t,R} - \langle v_i \rangle_{t,R} \langle v_j \rangle_{t,R})^2]\;;\\ 
\frac{\partial F_v}{\partial R_v}
&= \frac{\alpha_v}{n_v} \sum_{i=1}^{n_u} \sum_{j=1}^{n_v} \E[(\langle u_i v_j \rangle_{t,R} - \langle u_i \rangle_{t,R} \langle v_j \rangle_{t,R})^2]  \;.
\end{align*}
Both $\nicefrac{\partial F_v}{\partial R_u}$, $\nicefrac{\partial F_v}{\partial R_v}$ are clearly nonnegative. If $(\bu,\bv)$ are jointly distributed w.r.t.\ the posterior \eqref{posterior_H_t_R} then $\Vert \bu \Vert = \sqrt{\rho_u n_u}$ and $\Vert \bv \Vert = \sqrt{\rho_v n_v}$; hence $\nicefrac{\partial F_v}{\partial R_u} \leq 4 \alpha_v \rho_v^2 n_v$ and $\nicefrac{\partial F_v}{\partial R_v} \leq 4 \alpha_v \rho_u\rho_v n_u$. Thus, we have shown than $J_{F(t,\cdot)}(R)$ is uniformly bounded in $(t,R)$.

By Cauchy-Lipschitz theorem, for all $\epsilon=(\epsilon_u, \epsilon_v) \in [0, +\infty)^2$ there exists a unique solution $R(\cdot, \epsilon): [0,\delta] \to [0,+\infty)^2$ to the initial value problem $g' = F(t, g)$, $g(0)=\epsilon$.
Here $\delta \in [0,1]$ is such that $[0,\delta]$ is the maximal interval of existence of the solution.
Because $F$ has its image in $[0, \rho_u] \times [0, \rho_v]$, we have $R([0, \delta],\epsilon) \subseteq [\epsilon_u, \epsilon_u + \delta \rho_u] \times [\epsilon_v, \epsilon_v + \delta \rho_v]$, which means that $\delta = 1$ (the solution never leaves the domain of definition of $F$).

Each initial condition $\epsilon \in [0, +\infty)^2$ is tied to a unique solution $R(\cdot,\epsilon)$. This implies that the function $\epsilon \mapsto R(t,\epsilon)$ is injective. Its Jacobian determinant is given by Liouville's formula \cite{hartman1982ordinary}:
\begin{equation*}
\det J_{R(t,\cdot)}(\epsilon)
= \exp \int_0^t ds \, \bigg(\frac{\partial F_u}{\partial R_u} + \frac{\partial F_v}{\partial R_v}\bigg)\bigg\vert_{s,R(s,\epsilon)}\,
\end{equation*}
This Jacobian determinant is greater than, or equal to, one as both $\nicefrac{\partial F_u}{\partial R_u}$ and $\nicefrac{\partial F_v}{\partial R_v}$ are nonnegative.
For $\nicefrac{\partial F_v}{\partial R_v}$, this follows from our previous computations.
For $\nicefrac{\partial F_u}{\partial R_u}$, notice that $\nicefrac{\partial F_u}{\partial R_u} = \frac{\alpha_v \rho_u^2}{(1 + \alpha_v \rho_u F_v(t,R))^2}\nicefrac{\partial F_v}{\partial R_u}$ where $\nicefrac{\partial F_v}{\partial R_u} \geq 0$.
The fact that the Jacobian determinant is bounded away from $0$ uniformly in $\epsilon$ implies by the inverse function theorem that the injective function  $\epsilon \mapsto R(t,\epsilon)$ is a $\mathcal{C}^1$-diffeomorphism from $[0, +\infty)^2$ onto its image.
\end{IEEEproof}
\newpage
\section{Concentration of the overlap}\label{app:concentration_overlap}
Remember that the angular brackets $\langle - \rangle_{t,R}$ denote the expectation with respect to the posterior distribution \eqref{posterior_H_t_R}.
One central result in order to prove the lower bound on $\liminf_n f_n$ and the upper bound on $\limsup_n f_n$ is the concentration of the scalar overlap $Q_v$ around its expectation $\E \langle Q_v \rangle_{t,R}$, as long as we integrate over $R$ in a bounded subset of $[0,+\infty)^2$.
This corresponds to Proposition~\ref{prop:concentration_overlap} in the main text, that we repeat here for reader's convenience:
\begin{proposition*}[Concentration of the overlap around its expectation]
Let $M_u, M_v >0$.
For $n$ large enough, there exists a constant $M$ such that $\forall (a,b) \in (0,M_u)^2: a < \min\{1,b\}$, $\forall \delta \in (0,a)$, $\forall R_v \in [0,M_v]$, $\forall t \in [0,1]$:
\begin{equation*}
\int_{a}^{b} \E\,\big\langle \big(Q_v -\E\,\langle Q_v \rangle_{t,R}\,\big)^2\,\big\rangle_{t,R}\,dR_u
\leq M\bigg(\frac{1}{\delta^2 n} - \frac{\ln(a)}{n} + \frac{\delta}{a-\delta}\bigg)\;.
\end{equation*}
\end{proposition*}
The proof of Proposition~\ref{prop:concentration_overlap} is carried out mostly as in \cite{BarbierMacrisAdaptiveInterpolation}.
The main difference is that we don't need to assume that the marginals of the prior $P_v$ have a support bounded uniformly with $n$.
It will be enough that the norm of a vector distributed with respect to $P_v$ scales likes $\sqrt{n}$.
The concentration of the overlap around its expectation will follow from the concentration of the quantity:
\begin{equation}\label{def_L}
\cL
= \frac{1}{n}\sum_{j=1}^{n_v} \frac{\alpha_u}{2} v_j^2 - \alpha_u \, v_j V_j - \frac{1}{2}\sqrt{\frac{\alpha_u}{R_u}}\, v_j \overline{Z}_j \:.
\end{equation}
We first prove the following lemma that links the fluctuations of $Q_v$ to the fluctuations of $\cL$.
\begin{lemma}\label{lemma:computation_E<L>_and_others}
$\forall (t,R) \in [0,1] \times (0, +\infty)^2$, we have:
	\begin{align}
	\E\,\langle \cL \rangle_{t,R} &= -\frac{\alpha_u}{2} \frac{n_v}{n}\E\,\langle Q_v\rangle_{t,R}\;;\label{formula_E<L>}\\
 \E\,\langle (Q_v - \E\,\langle Q_v\rangle_{t,R})^2 \rangle_{t,R}
&\leq \frac{4}{\alpha_u^2} \bigg(\frac{n}{n_v}\bigg)^{\!\! 2}\, \E\,\langle (\cL - \E\,\langle \cL \rangle_{t,R})^2 \rangle_{t,R}\;.\label{upperbound_fluctuation_Q_v}
	\end{align}
\end{lemma}
\begin{IEEEproof}
Fix $(t,R) \in [0,1] \times (0,+\infty)^2$.
By definition \eqref{def_L} of $\cL$, we have:
\begin{align}
\E\,\langle \cL \rangle_{t,R}
	&= \frac{1}{n}\sum_{j=1}^{n_v} \frac{\alpha_u}{2} \E \langle v_j^2 \rangle_{t,R} - \alpha_u \, \E\big[ \langle v_j \rangle_{t,R} V_j \big] - \frac{1}{2}\sqrt{\frac{\alpha_u}{R_u}}\, \E\big[ \langle v_j \rangle_{t,R} \overline{Z}_j \big] \;;\label{def_E<L>}\\
\E\,\langle Q_v \cL \rangle_{t,R}
	&= \frac{1}{n}\sum_{j=1}^{n_v} \frac{\alpha_u}{2} \E \langle Q_v v_j^2 \rangle_{t,R} - \alpha_u \, \E\big[ \langle Q_v v_j \rangle_{t,R} V_j \big] - \frac{1}{2}\sqrt{\frac{\alpha_u}{R_u}}\, \E\big[ \langle Q_v v_j \rangle_{t,R} \overline{Z}_j \big]\,.\label{def_E<Q_v L>}
\end{align}
Integrating by parts with respect to the Gaussian random variable $\overline{Z}_j$, the last expectation on the right-hand side of each of \eqref{def_E<L>} and \eqref{def_E<Q_v L>} reads:
\begin{align}
\E\big[ \langle v_j \rangle_{t,R} \overline{Z}_j \big]
	&= \sqrt{\alpha_u R_u} \E\big[ \langle v_j^2 \rangle_{t,R} \big]
	-\sqrt{\alpha_u R_u} \E\big[ \langle v_j \rangle_{t,R}^2\big] \;;\label{stein_lemma_last_term_E<L>}\\
\E\big[ \langle Q_v v_j \rangle_{t,R} \overline{Z}_j \big]
	&= \sqrt{\alpha_u R_u} \E\big[ \langle Q_v v_j^2 \rangle_{t,R} \big]
	-\sqrt{\alpha_u R_u} \E\big[ \langle Q_v v_j \rangle_{t,R} \langle v_j \rangle_{t,R} \big] \;.\label{stein_lemma_last_term_E<Q_v L>}
\end{align}
Plugging \eqref{stein_lemma_last_term_E<L>} in \eqref{def_E<L>} yields:
\begin{equation*}
	\E\,\langle \cL \rangle_{t,R}
= \frac{\alpha_u}{n}\sum_{j=1}^{n_v} \frac{1}{2} \E\big[ \langle v_j \rangle_{t,R}^2\big] - \, \E\big[ \langle v_j \rangle_{t,R} V_j \big]
= -\frac{\alpha_u}{2} \frac{n_v}{n}\sum_{j=1}^{n_v} \frac{\E\big[ \langle v_j \rangle_{t,R} V_j \big]}{n_v}
= -\frac{\alpha_u}{2} \frac{n_v}{n} \E\,\langle Q_v \rangle_{t,R} \;,
\end{equation*}
where the second equality follows from Nishimori identity $\E[ \langle v_j \rangle_{t,R}^2] = \E[ \langle v_j \rangle_{t,R} V_j]$. This ends the proof of \eqref{formula_E<L>}. Plugging \eqref{stein_lemma_last_term_E<Q_v L>} in \eqref{def_E<Q_v L>}, it comes:
\begin{align}
\E\,\langle Q_v \cL \rangle_{t,R}
&= \frac{\alpha_u}{n}\sum_{j=1}^{n_v} \frac{1}{2}\E\big[ \langle Q_v v_j \rangle_{t,R} \langle v_j \rangle_{t,R} \big] - \E\big[ \langle Q_v v_j \rangle_{t,R} V_j \big]\nonumber\\
&= \frac{\alpha_u}{n}\sum_{j=1}^{n_v} \frac{1}{2}\E\big[ \langle Q_v \rangle_{t,R} \langle v_j V_j \rangle_{t,R} \big] - \E\big[ \langle Q_v v_j \rangle_{t,R} V_j \big]\nonumber\\
&= \alpha_u\frac{n_v}{n}\bigg( \frac{1}{2}\E\big[ \langle Q_v \rangle_{t,R}^2 \big] - \E\,\langle Q_v^2 \rangle_{t,R} \bigg) \;.\label{E<Q_v L>}
	\end{align}
The second equality follows once again from Nishimori identity
$$
\E\big[ \langle Q_v v_j \rangle_{t,R} \langle v_j \rangle_{t,R} \big]
= \frac{1}{n_v}\sum_{i=1}^{n_v}\E\big[ \langle v_i V_i v_j \rangle_{t,R} \langle v_j \rangle_{t,R} \big]
= \frac{1}{n_v}\sum_{i=1}^{n_v}\E\big[ V_i \langle  v_i \rangle_{t,R} V_j\langle v_j \rangle_{t,R} \big]
= \E\big[ \langle  Q_v \rangle_{t,R} \langle v_j V_j\rangle_{t,R} \big] \;.
$$
Combining \eqref{E<Q_v L>} and \eqref{formula_E<L>} yields:
\begin{align}
\E\,\langle Q_v (\cL - \E\,\langle \cL \rangle_{t,R}) \rangle_{t,R}
&= \E\,\langle Q_v \cL \rangle_{t,R} - \E\,\langle Q_v \rangle_{t,R} \E\,\langle \cL \rangle_{t,R}\nonumber\\
&= \frac{\alpha_u}{2}\frac{n_v}{n}\bigg(\E\big[ \langle Q_v \rangle_{t,R}^2 \big] - 2\E\,\langle Q_v^2 \rangle_{t,R} + \E[\langle Q_v\rangle_{t,R}]^2\bigg)\nonumber\\
&= -\frac{\alpha_u}{2}\frac{n_v}{n}\bigg(\E\,\big\langle(Q_v - \langle Q_v \rangle_{t,R})^2 \big\rangle_{t,R} + \E\,\big\langle (Q_v - \E\,\langle Q_v\rangle_{t,R})^2 \big\rangle_{t,R}\bigg) \;. \label{identity_E<Q(L - E<L>)>}
\end{align}
The identity \eqref{identity_E<Q(L - E<L>)>} directly implies (the second inequality below is an application of Cauchy-Schwarz inequality):
\begin{align}
\frac{\alpha_u}{2}\frac{n_v}{n} \E\,\big\langle (Q_v - \E\,\langle Q_v\rangle_{t,R})^2 \big\rangle_{t,R}
&\leq \big\vert \E\,\langle Q_v (\cL - \E\,\langle \cL \rangle_{t,R}) \rangle_{t,R} \big\vert\nonumber\\
&= \big\vert \E\,\langle (Q_v - \E\,\langle Q_v \rangle_{t,R})(\cL - \E\,\langle \cL \rangle_{t,R}) \rangle_{t,R} \big\vert\nonumber\\
&\leq \sqrt{\E\,\langle (Q_v - \E\,\langle Q_v \rangle_{t,R})^2 \rangle_{t,R} \cdot \E\,\langle (\cL - \E\,\langle \cL \rangle_{t,R})^2 \rangle_{t,R}}\;.
\end{align}
The upper bound \eqref{upperbound_fluctuation_Q_v} on the fluctuation of $Q_v$ follows simply from this last upper bound.
\end{IEEEproof}
\subsection{Concentration of $\cL$ around its expectation}
To prove concentration results on $\cL$, it will be useful to work with the free entropy $\frac{1}{n} \ln \cZ_{t,R}(\bY^{(t)},\widetilde{\bY}^{(t,R_v)}, \overline{\bY}^{(t,R_u)})$ where $\cZ_{t,R}(\bY^{(t)},\widetilde{\bY}^{(t,R_v)}, \overline{\bY}^{(t,R_u)})$ is the normalization factor of the Gibbs posterior distribution \eqref{posterior_H_t_R}.
In Appendix \ref{app:concentration_free_entropy}, we prove that this free entropy concentrates around its expectation when $n \to +\infty$. 
In order to shorten notations, we define:
\begin{equation}
F_n(t,R) \triangleq \frac{1}{n} \ln \cZ_{t,R}\big(\bY^{(t)},\widetilde{\bY}^{(t,R_v)}, \overline{\bY}^{(t,R_u)}\big)\;;\quad
f_n(t,R) \triangleq \frac{1}{n} \E\big[\ln \cZ_{t,R}\big(\bY^{(t)},\widetilde{\bY}^{(t,R_v)}, \overline{\bY}^{(t,R_u)}\big)\big] = \E\,F_n(t,R) \,.
\end{equation}
\begin{proposition}[Thermal fluctuations of $\cL$]\label{prop:concentration_L_on_<L>}
For $n$ large enough, we have for all positive real numbers $a < b$, $t \in [0,1]$ and $R_v \in [0,+\infty)$:
	\begin{equation}
	\int_a^b dR_u\,\E\,\big\langle \big(\cL - \langle \cL \rangle_{t,R}\big)^{2}\,\big\rangle_{t,R}
\leq \frac{\alpha_u \alpha_v \rho_v }{n} \bigg(\frac{\ln(b/a)}{2} + 1 \bigg)\,.
	\end{equation}
\end{proposition}
\begin{IEEEproof}
	Fix $(n,t) \in \mathbb{N}^* \times [0,1]$.
	Note that $\forall R \in (0, +\infty)^2$:
	\begin{equation}\label{1st_identity_derivative_fn_Rll'}
	\frac{\partial f_n}{\partial R_u}\bigg\vert_{t,R}
	= -\frac{1}{n}\E\Bigg[\Bigg\langle \frac{\partial \cH_{t,R}(\bx;\bY^{(t)},\widetilde{\bY}^{(t,R_v)}, \overline{\bY}^{(t,R_u)})}{\partial R_u} \Bigg\rangle_{\!\! t,R}\,\Bigg]
	=-\E\,\langle \cL \rangle_{t,R} \,.
	\end{equation}
	Further differentiating, we obtain:
	\begin{align}
		\frac{\partial^2 f_n}{\partial R_u^2}\bigg\vert_{t,R}
		&= \E\bigg[\bigg\langle \cL \,\frac{\partial \cH_{t,R}}{\partial R_u} \bigg\rangle_{\!\! t,R}\,\bigg]
		-\E\bigg[\langle \cL \rangle_{t,R}\,
		\bigg\langle \frac{\partial \cH_{t,R}}{\partial R_u} \bigg\rangle_{\!\! t,R}\,\bigg]
		-\E\,\bigg\langle \frac{\partial \cL }{\partial R_u} \bigg\rangle_{\!\! t,R}\\
		&=n\E\,\big\langle \big(\cL - \langle \cL \rangle_{t,R}\big)^{2}\,\big\rangle_{t,R}
		-\frac{1}{4 R_u }\sqrt{\frac{\alpha_u}{R_u}}\frac{\E\big[\langle \bv \rangle_{t,R}^T \overline{\bZ}\,\big]}{n} \;.\label{second_derivative_f_n_R}
	\end{align}
It follows directly from \eqref{second_derivative_f_n_R} that:
\begin{equation}\label{thermal_fluctuation_L}
\E\,\big\langle \big(\cL - \langle \cL \rangle_{t,R}\big)^{2}\,\big\rangle_{t,R}
	= \frac{1}{n}\frac{\partial^2 f_n}{\partial R_u^2}\bigg\vert_{t,R}
	+ \frac{1}{4 R_u }\sqrt{\frac{\alpha_u}{R_u}}\,\frac{\E\big[\langle \bv \rangle_{t,R}^\sT \overline{\bZ}\,\big]}{n^2}
\end{equation}
We start with upper bounding the integral over the second summand on the right-hand side of \eqref{thermal_fluctuation_L}.
Thanks to an integration by parts with respect to $\overline{Z}_j$, $j \in \{1,\dots,n_u\}$, it comes:
\begin{equation}
\frac{1}{4 R_u }\sqrt{\frac{\alpha_u}{R_u}}\,\frac{\E\big[\langle \bv \rangle_{t,R} \overline{\bZ}\,\big]}{n^2}
= \frac{\alpha_u}{4 R_u } \frac{\E\, \langle \Vert \bv \Vert^2 \rangle_{t,R} - \E\,\Vert \langle \bv  \rangle_{t,R}\Vert^2}{n^2}
\leq \frac{\alpha_u \rho_v}{4 R_u } \frac{n_v}{n^2}\:.
\end{equation}
Therefore:
	\begin{equation}\label{upperbound_int_dL/dR}
\int_a^b \frac{dR_u}{4 R_u }\sqrt{\frac{\alpha_u}{R_u}}\,\frac{\E\big[\langle \bv \rangle_{t,R}^T \overline{\bZ}\,\big]}{n^2}
	\leq \frac{n_v}{n^2} \frac{\alpha_u \rho_v \ln(b/a)}{4}\,.
	\end{equation}
	It remains to upper bound $\int_a^b \frac{dR_u}{n}\frac{\partial^2 f_n}{\partial R_u^2}\big\vert_{t,R}
			= \frac{1}{n}\frac{\partial f_n}{\partial R_u}\big\vert_{t,R_u=b, R_v} - \frac{1}{n}\frac{\partial f_n}{\partial R_u}\big\vert_{t,R_u=a, R_v}$.
	Note that $\forall R \in [0,+\infty)^2$:
	\begin{equation}\label{2nd_identity_derivative_fn_Rll'}
		\frac{\partial f_n}{\partial R_u}\bigg\vert_{t,R}
		= - \E\,\langle \cL \rangle_{t,R}
		= \frac{\alpha_u}{2}\frac{n_v}{n} \E\,\langle Q_v \rangle_{t,R}
		= \frac{\alpha_u}{2}\frac{n_v}{n} \frac{\E\big[\Vert\langle \bv \rangle_{t,R}\Vert^2\big]}{n_v}\:,
	\end{equation}
	where the first equality follows from \eqref{1st_identity_derivative_fn_Rll'}, the second one from Lemma \ref{lemma:computation_E<L>_and_others} (this is the identity \eqref{formula_E<L>}) and the third one from Nishimori identity: $\E\,\langle Q_v \rangle_{t,R} =  \nicefrac{\E[\langle \bv \rangle_{t,R}^{\sT}\bV]}{n_v} = \nicefrac{\E[\Vert\langle \bv \rangle_{t,R}\Vert^2]}{n_v}$.
	Making use of \eqref{2nd_identity_derivative_fn_Rll'} and Jensen's inequality (for the upper bound), it is easy to see that $\forall R \in [0,+\infty)^2$:
	\begin{equation}\label{upperbound_derivative_fn_Rll'}
	0 \leq \frac{\partial f_n}{\partial R_u}\bigg\vert_{t,R}
	\leq \frac{\alpha_u}{2}\frac{n_v}{n} \frac{\E\,\langle \Vert\bv\Vert^2 \rangle_{t,R}}{n_v}
	= \frac{\alpha_u \rho_v }{2}\frac{n_v}{n} \;.
	\end{equation}
	Combining both \eqref{upperbound_int_dL/dR} and \eqref{upperbound_derivative_fn_Rll'}, we finally get:
	\begin{equation}\label{final_upperbound_int_secondDerivative_fn_Rll'}
	\int_a^b dR_u\, \E\,\big\langle \big(\cL - \langle \cL \rangle_{t,R}\big)^{2}\,\big\rangle_{t,R}
	\leq \frac{1}{n} \frac{n_v}{n} \frac{\alpha_u \rho_v}{2}\bigg(\frac{\ln(b/a)}{2} + 1 \bigg)\;.
	\end{equation}
\end{IEEEproof}

\begin{proposition}[Quenched fluctuations of $\cL$]\label{prop:concentration_<L>_on_E<L>}
Let $M_u, M_v >0$.
For $n$ large enough, there exists a constant $M$ such that $\forall (a,b) \in (0,M_u)^2: a < \min\{1,b\}$, $\forall \delta \in (0,a)$, $\forall R_v \in [0,M_v]$, $\forall t \in [0,1]$:
\begin{equation}
	\int_{a}^{b} dR_u\,\E\,\big\langle \big(\langle \cL \rangle_{t,R}-\E\,\langle \cL \rangle_{t,R}\,\big)^2\,\big\rangle_{t,R}
	\leq M\bigg(\frac{1}{\delta^2 n} - \frac{\ln(a)}{n} + \frac{\delta}{a-\delta} \bigg)\,.
\end{equation}
\end{proposition}
\begin{IEEEproof}
	Fix $(n,t) \in \mathbb{N}^* \times [0,1]$. For all $R \in (0,+\infty)^2$, we have:
	\begin{align}
		\frac{\partial F_n}{\partial R_u}\bigg\vert_{t,R}
		&=-\langle \cL \rangle_{t,R} \;;\\
		\frac{\partial^2 F_n}{\partial R_u^2}\bigg\vert_{t,R}
		&=n \big\langle \big(\cL - \langle \cL\rangle_{t,R}\big)^{2}\,\big\rangle_{t,R}
	-\frac{1}{4 R_u }\sqrt{\frac{\alpha_u}{R_u}}\frac{\langle \bv \rangle_{t,R}^T\, \overline{\bZ}}{n}\;;\label{quenched:2ndDeriv_Fn}\\
		\frac{\partial f_n}{\partial R_u}\bigg\vert_{t,R}
		&=-\E\,\langle \cL \rangle_{t,R}\;;\\
		\frac{\partial^2 f_n}{\partial R_u^2}\bigg\vert_{t,R}
		&=n \E\,\big\langle \big(\cL - \langle \cL\rangle_{t,R}\big)^{2}\,\big\rangle_{t,R}
	-\frac{1}{4 R_u }\sqrt{\frac{\alpha_u}{R_u}}\frac{\E\big[\langle \bv \rangle_{t,R}^T \overline{\bZ}\,\big]}{n}\,.
	\end{align}
	The second term on the right-hand side of \eqref{quenched:2ndDeriv_Fn} can be upper bounded with Cauchy-Schwarz inequality:
	\begin{equation}\label{upperbound_2ndterm_2ndDeriv_Fn}
	\Bigg\vert\frac{1}{4 R_u }\sqrt{\frac{\alpha_u}{R_u}}\frac{\langle \bv \rangle_{t,R}^T\, \overline{\bZ}}{n}\Bigg\vert
	\leq \frac{1}{4 R_u }\sqrt{\frac{\alpha_u}{R_u}}\frac{\Vert \langle \bv \rangle_{t,R} \Vert\, \Vert \overline{\bZ} \Vert}{n}
	\leq \frac{1}{4 R_u }\sqrt{\frac{\alpha_u}{R_u}}\frac{\langle \Vert \bv \Vert \rangle_{t,R}\, \Vert \overline{\bZ} \Vert}{n}
	\leq \frac{1}{4 R_u }\sqrt{\frac{\alpha_u \rho_v}{R_u} \frac{n_v}{n}}\frac{\Vert \overline{\bZ} \Vert}{\sqrt{n}}\,.
	\end{equation}
	We now define for all $R_u \in (0, +\infty)$:
	\begin{align}
		F(R_u) &\triangleq F_n(t,(R_u,R_v)) - \sqrt{\alpha_u \rho_v R_u \frac{n_v}{n}}\frac{\Vert \overline{\bZ} \Vert}{\sqrt{n}}\:;\\
		f(R_u) &\triangleq f_n(t,(R_u,R_v)) - \sqrt{\alpha_u \rho_v R_u \frac{n_v}{n}}\frac{\E\,\Vert \overline{\bZ} \Vert}{\sqrt{n}}\,.
	\end{align}
	$F$ is convex on $(0, +\infty)$ as it is twice differentiable with a nonnegative second derivative by \eqref{quenched:2ndDeriv_Fn} and \eqref{upperbound_2ndterm_2ndDeriv_Fn}.
	The same holds for $f$.
	We will apply the following standard result to these two convex functions (we refer to \cite{BarbierMacrisAdaptiveInterpolation} for the proof):
	\begin{lemma}[An upper bound for differentiable convex functions]\label{lemma:diff_convex_functions}
		Let $g$ and $G$ be two differentiable convex functions defined on an interval $I\subseteq \mathbb{R}$.
		Let $r \in I$ and $\delta > 0$ such that $r \pm \delta \in I$. Then
		\begin{equation}
		\vert G'(r) - g'(r) \vert
		\leq C_{\delta}(r) + \frac{1}{\delta}\sum_{u \in \{-\delta,0,\delta\}} \vert G(r+u) - g(r+u) \vert \,,
		\end{equation}
		where $C_{\delta}(r) = g'(r+\delta) - g'(r-\delta) \geq 0$.
	\end{lemma}
	\noindent For all $R_u \in (0, +\infty)$, we have:
	\begin{align}
		F(R_u) - f(R_u)
		= F_n(t,(R_u, R_v)) - f_n(t,(R_u, R_v))
		- \sqrt{\alpha_u \rho_v R_u \frac{n_v}{n}}\frac{\Vert \overline{\bZ} \Vert - \E\,\Vert \overline{\bZ} \Vert}{\sqrt{n}} \:;\label{F_minus_f}\\
		F'(R_u) - f'(R_u)
		= -\Big(\langle \cL \rangle_{t,R} - \E\,\langle \cL\rangle_{t,R}\Big)
		- \frac{1}{2}\sqrt{\frac{\alpha_u \rho_v}{R_u} \frac{n_v}{n}}\frac{\Vert \overline{\bZ} \Vert - \E\,\Vert \overline{\bZ} \Vert}{\sqrt{n}} \;. \label{F'_minus_f'}
	\end{align}
	Let $C_{\delta}(r) = f'(r+\delta) - f'(r-\delta)$, which is nonnegative by convexity of $f$.
	It follows from Lemma \ref{lemma:diff_convex_functions} and the two identities \eqref{F_minus_f} and \eqref{F'_minus_f'} that $\forall R_u \in (0, +\infty)$, $\forall \delta \in (0,R_u)$:
	\begin{align*}
		\big\vert \langle \cL \rangle_{t,R} - \E\,\langle\cL \rangle_{t,R}\big\vert
		&\leq
		\frac{1}{2}\sqrt{\frac{\alpha_u \rho_v}{R_u} \frac{n_v}{n}}\frac{\big\vert \Vert \overline{\bZ} \Vert - \E\,\Vert \overline{\bZ} \Vert \big\vert}{\sqrt{n}}
		+ C_{\delta}(R_u)
		+ \frac{1}{\delta}\sum_{x \in \{-\delta,0,\delta\}} \vert F(R_u+x) - f(R_u+x)\vert\\
		&\leq
		\sqrt{\alpha_u \rho_v \frac{n_v}{n}}
		\bigg(\frac{1}{2\sqrt{R_u}} + 3\sqrt{R_u}\bigg)\frac{\big\vert \Vert \overline{\bZ} \Vert - \E\,\Vert \overline{\bZ} \Vert \big\vert}{\sqrt{n}}
		+ C_{\delta}(R_u)\\
		&\qquad\qquad\qquad\qquad\qquad\qquad\quad\, + \frac{1}{\delta}\sum_{x \in \{-\delta,0,\delta\}} \vert F_n(t,(R_u + x, R_v)) - f_n(t,(R_u + x, R_v))\vert \;.
	\end{align*}
	Thanks to the inequality $(\sum_{i=1}^{m} v_i)^2 \leq m \sum_{i=1}^{m} v_i^2$, this directly implies $\forall R_u \in (0, +\infty)$, $\forall \delta \in (0,R_u)$:
	\begin{multline}\label{upperbound_variance_thermal_L}
		\E\big[\big( \langle \cL \rangle_{t,R} - \E\,\langle\cL \rangle_{t,R}\big)^{2}\,\big]
		\leq 		5 \alpha_u \rho_v \frac{n_v}{n}
		\bigg(\frac{1}{4 R_u} + 3 + 9 R_u\bigg) \, \frac{\Var \Vert \overline{\bZ} \Vert}{n}
		+ 5C_{\delta}(R_u)^2\\
		+ \frac{5}{\delta^2}\sum_{x \in \{-\delta,0,\delta\}} \E\big[\big(F_n(t,(R_u + x, R_v)) - f_n(t,(R_u + x, R_v))\big)^2\big]\,.
	\end{multline}
	The next step is to bound the integral of the three summands on the right-hand side of \eqref{upperbound_variance_thermal_L}.
	By \cite[Theorem 3.1.1]{vershynin_2018}, there exists $C_1$ such that $\Var\,\Vert\overline{\bZ}\Vert \leq C_1$ independently of the dimension $n_v$. Then:
	\begin{equation}\label{upperbound_1st_summand}
	\int_{a}^{b} dR_u \,5 \alpha_u \rho_v \frac{n_v}{n}
	\bigg(\frac{1}{4 R_u} + 3 + 9 R_u\bigg) \, \frac{\Var \Vert \overline{\bZ} \Vert}{n}
	\leq 5 \alpha_u \rho_v \frac{n_v}{n}
	\bigg(\frac{\ln(b/a)}{4} + 3b + \frac{9}{2} b^2 \bigg) \, \frac{C_1}{n}
	\,.
	\end{equation}
	Note that $C_{\delta}(R_u) = \vert C_{\delta}(R_u)\vert \leq \vert f'(R_u+\delta) \vert + \vert f'(R_u-\delta) \vert$. For all $R_u \in (0, +\infty)$, we have:
	\begin{equation}\label{upperbound_f'}
		\vert f'(R_u) \vert
		\leq
		\big\vert \E\,\langle \cL \rangle_{t,R} \big\vert
		+ \frac{1}{2}\sqrt{\frac{\alpha_u \rho_v}{R_u} \frac{n_v}{n}}\frac{\E\,\Vert \overline{\bZ} \Vert}{\sqrt{n}}
		\leq \frac{n_v}{n} \frac{\sqrt{\alpha_u \rho_v}}{2} \bigg(\sqrt{\alpha_u \rho_v} + \frac{1}{\sqrt{R_u}}\bigg)\,,
	\end{equation}
	The second inequality in \eqref{upperbound_f'} follows from
	the upper bounds
	$\vert \E\,\langle \cL \rangle_{t,R} \vert \leq \nicefrac{\alpha_u \rho_v n_v}{2n}$ (see \eqref{upperbound_derivative_fn_Rll'}) and
	$\E \Vert\overline{\bZ}\Vert \leq \E[\Vert \overline{\bZ} \Vert^2]^{\nicefrac{1}{2}} = \sqrt{n_v}$.
	Thus, for the second summand, we obtain $\forall \delta \in (0,a)$:
	\begin{align}\label{upperbound_2nd_summand}
		\int_a^b dR_u \, C_{\delta}(R_u)^2
		&\leq \frac{n_v}{n} \frac{\sqrt{\alpha_u\rho_v}}{2} \bigg(\sqrt{\alpha_u \rho_v} + \frac{1}{\sqrt{a-\delta}} \,\bigg) \int_a^b dR_u \, C_{\delta}(R_u)\nonumber\\
		&= \frac{n_v}{n} \frac{\sqrt{\alpha_u\rho_v}}{2} \bigg(\sqrt{\alpha_u \rho_v} + \frac{1}{\sqrt{a-\delta}} \,\bigg)
		\big[\big(f(b+\delta) - f(b-\delta)\big) - \big(f(a + \delta) - f(a -\delta)\big)\big]\nonumber\\
		&\leq \delta \bigg(\frac{n_v}{n}\bigg)^{\!\! 2} \alpha_u \rho_v \bigg(\sqrt{\alpha_u \rho_v} + \frac{1}{\sqrt{a-\delta}} \,\bigg)^{\! 2}\:.
	\end{align}
	The last inequality is a simple application of the mean value theorem.
	We finally turn to the third summand.
	By Proposition~\ref{prop:concentration_free_entropy} in Appendix~\ref{app:concentration_free_entropy}, there exists a positive constant $C_2$ depending only on $a$, $b$ and $M_v$ such that $\forall t \in [0,1]$, $\forall (R_u, R_v) \in (0,b+a) \times (0, M_v)$:
	\begin{equation}\label{upperbound_variance_free_entropy}
	\E\big[\big(F_n(t,R) - f_n(t,R)\big)^2 \,\big] \leq \frac{C_2}{n}\,.
	\end{equation}
	Using \eqref{upperbound_variance_free_entropy}, we see that the third summand satisfies $\forall \delta \in (0,a)$:
	\begin{equation}\label{upperbound_3rd_summand}
		\int_{a}^{b} \! dR_u \,
		\frac{5}{\delta^2}\sum_{x \in \{-\delta,0,\delta\}} \E\big[\big(F_n(t,(R_u+x, R_v)) - f_n(t,(R_u + x, R_v))\big)^2 \,\big]
		\leq \frac{15C_2}{\delta^2 n} b \:.
	\end{equation}
	To end the proof it remains to integrate \eqref{upperbound_variance_thermal_L} over $R_u \in [a,b]$ and use the three upper bounds \eqref{upperbound_1st_summand}, \eqref{upperbound_2nd_summand} and \eqref{upperbound_3rd_summand}.
\end{IEEEproof}
\subsection{Concentration of $Q_v$ around its expectation: proof of Proposition~\ref{prop:concentration_overlap}}
Using the upper bound \eqref{upperbound_fluctuation_Q_v} and Cauchy-Schwarz inequality, it directly comes:
\begin{equation}
\int_{a}^{b} \E\,\big\langle \big( Q_v -\E\,\langle Q_v \rangle_{t,R}\,\big)^2\,\big\rangle_{t,R}\,dR_u
\leq \frac{4}{\alpha_u^2} \bigg(\frac{n}{n_v}\bigg)^{\!\! 2} \int_a^b \E\,\langle (\cL - \E\,\langle \cL \rangle_{t,R})^2 \rangle_{t,R}\,dR_u \;.
\end{equation}
We then use the concentration results for $\cL$, that is, Propositions~\ref{prop:concentration_L_on_<L>} and~\ref{prop:concentration_<L>_on_E<L>}, to upper bound
$$
\int_a^b \E\,\langle (\cL - \E\,\langle \cL \rangle_{t,R})^2 \rangle_{t,R}\,dR_u
= \int_a^b \E\,\langle (\cL - \langle \cL \rangle_{t,R})^2 \rangle_{t,R}\,dR_u + \int_a^b \E[(\langle \cL \rangle_{t,R} - \E\,\langle \cL \rangle_{t,R})^2\,]\,dR_u
$$
and prove Proposition~\ref{prop:concentration_overlap}.

\newpage
\section{Concentration of the free entropy}\label{app:concentration_free_entropy}
Consider the inference problem~\eqref{interpolation_model_R}.
Once the observations $\bY^{(t)}$, $\widetilde{\bY}^{(t,R_v)}$ and $\overline{\bY}^{(t,R_u)}$ have been replaced by their definitions, the associated Hamiltonian reads:
\begin{multline}
\cH_{t,R}(\bu, \bv ; \bU, \bV, \bZ, \widetilde{\bZ}, \overline{\bZ})
\triangleq
\sum_{i=1}^{n_u}\sum_{j=1}^{n_v}
\frac{(1-t)}{2n} u_i^2 v_j^2 - \frac{1-t}{n}\, u_i v_j U_i V_j - \sqrt{\frac{1-t}{n}}u_i v_j Z_{ij}\\
+
\sum_{i=1}^{n_u} \frac{\alpha_v R_v}{2} u_i^2 - \alpha_v R_v \, u_i U_i- \sqrt{\alpha_v R_v}\, u_i \widetilde{Z}_i
+
\sum_{j=1}^{n_v} \frac{\alpha_u R_u}{2} v_j^2 - \alpha_u R_u v_j V_j - \sqrt{\alpha_u R_u}\, v_j \overline{Z}_j\:.
\end{multline}
In this section, we show that the free entropy
\begin{equation}
\frac{1}{n} \ln \cZ_{t,R}\big(\bY^{(t)},\widetilde{\bY}^{(t,R_v)}, \overline{\bY}^{(t,R_u)}\big)
= \frac{1}{n} \ln\Bigg( \int dP_u(\bu) dP_v(\bv)\: e^{-\cH_{t,R}(\bu, \bv ; \bU, \bV, \bZ, \widetilde{\bZ}, \overline{\bZ})}\Bigg)
\end{equation}
concentrates around its expectation.
We will sometimes write $\frac{1}{n} \ln \cZ_{t,R}$, omitting the arguments, to shorten notations.
\begin{proposition}[Concentration of the free entropy]\label{prop:concentration_free_entropy}
For any positive number $M$, there exists a positive constant $C$ such that for any $R \in [0,+\infty)^2$ whose Euclidian norm is bounded by $M$ we have: 
\begin{equation}\label{bound_variance_free_entropy}
\E \Bigg[\Bigg(\frac{ \ln \cZ_{t,R}}{n}
	- \E\bigg[\frac{\ln \cZ_{t,R}}{n} \bigg]
	\Bigg)^{\!\! 2}\:\Bigg]
	\leq \frac{C}{n} \:.
\end{equation}
\end{proposition}
\begin{IEEEproof}
We drop the subscripts to the angular brackets $\langle - \rangle_{t,R}$ to lighten notations.
First, we show that the free entropy concentrates on its conditional expectation given $\bV$, $\bZ$, $\widetilde{\bZ}$, $\overline{\bZ}$.
Thus, $\nicefrac{\ln  \cZ_{t,R}}{n}$ is seen as a function of $\nicefrac{\bU}{\sqrt{\rho_u n_u}}$ and we work conditionally to $\bV$, $\bZ$, $\widetilde{\bZ}$, $\overline{\bZ}$: $g(\nicefrac{\bU}{\sqrt{\rho_u n_u}}) \equiv \nicefrac{\ln  \cZ_{t,R}}{n}$.
We normalize by $\sqrt{\rho_u n_u}$ because $\nicefrac{\bU}{\sqrt{\rho_u n_u}}$ is uniformly distributed on the $(n_u-1)$-sphere of radius $1$, and we want to apply L\'{e}vy's lemma on the concentration of uniform measure on the sphere.
A statement of this lemma, which we reproduce here for the reader's convenience, can be found in \cite[Corollary 5.4]{haar_meckes} along with a proof.
\begin{lemma}[L\'{e}vy's lemma]\label{lemma:concentration_vector_sphere}
	Let $\mathcal{S}^{n-1}$ the $(n-1)$-sphere of radius $1$.
	Let $f: \mathcal{S}^{n-1} \to \R$ be Lipschitz with Lipschitz constant $L$, and let $\bX$ be a uniform random vector in $\mathcal{S}^{n-1}$.
	Then
	$$
	\mathbb{P}(\vert f(\bX) - \E f(\bX)\vert \geq Lt) \leq \exp(\pi - n t^2/4) \;.
	$$
\end{lemma}
By Jensen's inequality, we have:
\begin{multline}\label{Jensen_difference_freee_entropy}
\frac{1}{n} \langle \cH_{t,R}(\bu, \bv ; \widetilde{\bU}, \bV, \bZ, \widetilde{\bZ}, \overline{\bZ})
- \cH_{t,R}(\bu, \bv ; \bU, \bV, \bZ, \widetilde{\bZ}, \overline{\bZ})\rangle_{\widetilde{\bU}}\\
\leq g\bigg(\frac{\bU}{\sqrt{\rho_u n_u}}\bigg) - g\bigg(\frac{\widetilde{\bU}}{\sqrt{\rho_u n_u}}\bigg)\\
\leq \frac{1}{n} \langle \cH_{t,R}(\bu, \bv ; \widetilde{\bU}, \bV, \bZ, \widetilde{\bZ}, \overline{\bZ})
- \cH_{t,R}(\bu, \bv ; \bU, \bV, \bZ, \widetilde{\bZ}, \overline{\bZ})\rangle_{\bU}
\end{multline}
The subscript $\widetilde{\bU}$ (resp.\ $\bU$) of the Gibbs bracket notation on the left-hand side (resp.\ right-hand side) of \eqref{Jensen_difference_freee_entropy} specifies that $(\bu, \bv)$ is distributed according to $dP_u(\bu) dP_v(\bv)\: e^{-\cH_{t,R}(\bu, \bv ; \widetilde{\bU}, \bV, \bZ, \tilde{\bZ}, \overline{\bZ})}$ (resp.\ $dP_u(\bu) dP_v(\bv)\: e^{-\cH_{t,R}(\bu, \bv ; \bU, \bV, \bZ, \tilde{\bZ}, \overline{\bZ})}$).
Note that
\begin{align*}
\big\vert \cH_{t,R}(\bu, \bv ; \widetilde{\bU}, \bV, \bZ, \widetilde{\bZ}, \overline{\bZ})
- \cH_{t,R}(\bu, \bv ; \bU, \bV, \bZ, \widetilde{\bZ}, \overline{\bZ})\big\vert
&=
\bigg\vert \frac{1-t}{n} \sum_{i=1}^{n_u} \sum_{j=1}^{n_v} v_j V_j u_i (U_i - \widetilde{U}_i)
+
\alpha_v R_v \sum_{i=1}^{n_u} u_i (U_i - \widetilde{U}_i) \bigg\vert\\
&=
\bigg\vert \frac{1-t}{n} \bv^{\sT} \bV
+
\alpha_v R_v\bigg\vert \cdot \big\vert \bu^{\sT} (\bU - \widetilde{\bU}) \big\vert \\
&\leq \bigg( \frac{n_v}{n} \rho_v
+
\alpha_v R_v \bigg) \rho_u n_u \bigg\Vert \frac{\bU}{\sqrt{\rho_u n_u}} - \frac{\widetilde{\bU}}{\sqrt{\rho_u n_u}} \bigg\Vert \;.
\end{align*}
Combining this last inequality with \eqref{Jensen_difference_freee_entropy} yields:
\begin{equation*}
\bigg\vert g\bigg(\frac{\bU}{\sqrt{\rho_u n_u}}\bigg) - g\bigg(\frac{\widetilde{\bU}}{\sqrt{\rho_u n_u}}\bigg)
\bigg\vert
\leq \rho_u \frac{n_u}{n} \bigg( \rho_v \frac{n_v}{n}
+
\alpha_v R_v \bigg) \bigg\Vert \frac{\bU}{\sqrt{\rho_u n_u}} - \frac{\widetilde{\bU}}{\sqrt{\rho_u n_u}} \bigg\Vert \;,
\end{equation*}
i.e., $g$ is Lipschitz with Lipschitz constant $L = \rho_u \frac{n_u}{n} \big( \rho_v \frac{n_v}{n}
+ \alpha_v R_v \big)$. Lemma~\ref{lemma:concentration_vector_sphere} directly implies:
\begin{equation}\label{bound_variance_Levy_1}
\E \Bigg[\Bigg(\frac{\ln  \cZ_{t,R}}{n} - \E\bigg[\frac{ \ln \cZ_{t,R}}{n} \bigg\vert \bV, \bZ, \widetilde{\bZ}, \overline{\bZ}\bigg]\Bigg)^{\!\! 2}\,\Bigg]
\leq \frac{4 L^2 e^{\pi}}{n_u}
= \frac{C_1}{n}
\end{equation}
with $C_1 = 4 e^{\pi} \rho_u^2 \frac{n_u}{n} \big( \rho_v \frac{n_v}{n}
+
\alpha_v R_v \big)^2$.\\

We can show in a similar way that the conditional expectation of the free entropy given $\bV$, $\bZ$, $\widetilde{\bZ}$, $\overline{\bZ}$ concentrates on its conditional expectation given $\bZ$, $\widetilde{\bZ}$, $\overline{\bZ}$, that is:
\begin{equation}\label{bound_variance_Levy_2}
\E \Bigg[\Bigg(\E\bigg[\frac{ \ln \cZ_{t,R}}{n} \bigg\vert \bV, \bZ, \widetilde{\bZ}, \overline{\bZ}\bigg] - \E\bigg[\frac{ \ln \cZ_{t,R}}{n} \bigg\vert \bZ, \widetilde{\bZ}, \overline{\bZ}\bigg]\Bigg)^{\!\! 2}\,\Bigg]
\leq \frac{C_2}{n}
\end{equation}
with $C_2 = 4 e^{\pi} \rho_v^2 \frac{n_v}{n} \big( \rho_u \frac{n_u}{n}
+
\alpha_u R_u \big)^2$.\\

Finally, we show that the conditional expectation of the free entropy given $\bZ$, $\widetilde{\bZ}$, $\overline{\bZ}$ concentrates on its expectation.
$\nicefrac{\ln  \cZ_{t,R}}{n}$ is seen as a function of the Gaussian noises $\bZ$, $\widetilde{\bZ}$, $\overline{\bZ}$: let $g(\bZ,\widetilde{\bZ}, \overline{\bZ}) \equiv \E[\nicefrac{\ln  \cZ_{t,R}}{n} \vert \bZ, \widetilde{\bZ}, \overline{\bZ}]$.
By the Gaussian-Poincar\'{e} inequality (see \cite[Theorem 3.20]{boucheron_concentration}), we have:
		\begin{equation}\label{gaussian_poincare}
		\E \Bigg[\Bigg(\E\bigg[\frac{\ln  \cZ_{t,R}}{n} \bigg\vert \bZ, \widetilde{\bZ}, \overline{\bZ}\bigg] - \E\bigg[\frac{ \ln \cZ_{t,R}}{n}\bigg]\Bigg)^{\!\! 2}\,\Bigg]
		\leq \E\,\big\Vert \nabla g(\bZ,\widetilde{\bZ}, \overline{\bZ}) \big\Vert^2 \;.
		\end{equation}
		The squared norm of the gradient of $g$ reads
		$\Vert \nabla g\Vert^2 =
		\sum_{i, j} \vert\nicefrac{\partial g}{\partial Z_{i,j}}\vert^2
		+ \sum_{i} \vert\nicefrac{\partial g}{\partial \widetilde{Z}_{i}}\vert^2
		+ \sum_{j} \vert\nicefrac{\partial g}{\partial \overline{Z}_{j}}\vert^2$.
		Each of these partial derivatives takes the form 
		$\nicefrac{\partial g}{\partial x} = -n^{-1} \big\langle \nicefrac{\partial \mathcal{H}_{t,R}}{\partial x} \big\rangle$. More precisely:
		\begin{equation*}
		\bigg\vert\frac{\partial g}{\partial Z_{ij}}\bigg\vert
		= n^{-1}\bigg\vert \sqrt{\frac{1-t}{n}} \langle  u_i v_j \rangle \bigg\vert\quad ; \quad
		\bigg\vert \frac{\partial g}{\partial \widetilde{Z}_{i}} \bigg\vert
		= n^{-1} \big\vert \sqrt{\alpha_v R_v} \, \langle u_i \rangle \big\vert\quad ; \quad
		\bigg\vert \frac{\partial g}{\partial \overline{Z}_{j}} \bigg\vert
		= n^{-1} \big\vert \sqrt{\alpha_u R_u} \, \langle v_j \rangle \big\vert \,.
		\end{equation*}
		On one hand, we have
		\begin{equation}\label{first_upperbound_GP}
		\sum_{i=1}^{n_u} \sum_{j=1}^{n_v} \E\,\bigg\vert\frac{\partial g}{\partial Z_{ij}}\bigg\vert^2
		\leq \frac{1}{n^{3}}\sum_{i=1}^{n_u} \sum_{j=1}^{n_v}
		\E[\langle u_i v_j \rangle^{2}]
		\leq \frac{1}{n^{3}}\sum_{i=1}^{n_u} \sum_{j=1}^{n_v}
		\E[\langle u_i^2 v_j^2 \rangle]
		= \frac{n_u}{n} \frac{n_v}{n} \frac{\rho_u \rho_v}{n}\:,
		\end{equation}
		where the second inequality follows from Jensen's inequality.
		On the other hand, we have
		\begin{equation}\label{second_upperbound_GP}
		\sum_{i=1}^{n_u} \E\,\bigg\vert\frac{\partial g}{\partial \widetilde{Z}_{i}}\bigg\vert^2
		\leq \frac{\alpha_v R_v}{n^2} \sum_{i=1}^{n_u} \E[\langle u_i^2 \rangle]
		= \frac{n_u}{n}\frac{\alpha_v \rho_u R_v}{n} \quad ; \quad
		\sum_{j=1}^{n_v} \E\,\bigg\vert\frac{\partial g}{\partial \overline{Z}_{j}}\bigg\vert^2
		\leq \frac{\alpha_u R_u}{n^2} \sum_{j=1}^{n_v} \E[\langle v_j^2 \rangle]
		= \frac{n_v}{n}\frac{\alpha_u \rho_v R_u}{n}\;;
		\end{equation}
		where in both cases the first inequality follows from Jensen's inequality once again.
		Plugging \eqref{first_upperbound_GP} and \eqref{second_upperbound_GP} in \eqref{gaussian_poincare} yields:
		\begin{equation}\label{bound_variance_GP}
		\E \Bigg[\Bigg(\E\bigg[\frac{\ln  \cZ_{t,R}}{n} \bigg\vert \bZ, \widetilde{\bZ}, \overline{\bZ}\bigg] - \E\bigg[\frac{ \ln \cZ_{t,R}}{n}\bigg]\Bigg)^{\!\! 2}\,\Bigg]
		\leq \frac{C_3}{n}\,,
		\end{equation}
		with $C_3 = \frac{n_u}{n} \frac{n_v}{n} \rho_u \rho_v + \frac{n_u}{n} \alpha_v \rho_u R_v + \frac{n_v}{n} \alpha_u \rho_v R_u$.\\
Note that
$$
C_1 + C_2 + C_3 \xrightarrow[n \to +\infty]{} C \triangleq \alpha_u \alpha_v (4 e^\pi \alpha_v \rho_u^2 (\rho_v + R_v)^2
+ 4 e^\pi \alpha_u \rho_v^2 (\rho_u + R_u)^2 + \rho_u \rho_v + \rho_u R_v + \rho_v R_v) \;.
$$
This limit combined with the inequalities \eqref{bound_variance_Levy_1},  \eqref{bound_variance_Levy_2}, and \eqref{bound_variance_GP} ends the proof of \eqref{bound_variance_free_entropy}.
\end{IEEEproof}
\newpage
\section{Proof of Theorem~\ref{theorem:MMSE}: formula for the asymptotic matrix-MMSE}\label{app:proof_mmse}
In the whole appendix we suppose that the positive hyperparameters $\alpha_u$, $\alpha_v$, $\rho_u$ and $\rho_v$ are all fixed.
Then, we define $\forall (m_u,m_v, \lambda) \in [0,\rho_u] \times [0,\rho_v] \times (0,+\infty)$:
\begin{equation*}
i(m_u, m_v, \lambda)
\triangleq i_{\scriptscriptstyle \Theta}(m_u,m_v)
= \frac{\lambda \alpha_u \alpha_v}{2}(\rho_u-m_u)(\rho_v-m_v)
+ \alpha_u \frac{\ln(1+\lambda\alpha_v \rho_u m_v)}{2}
+ \alpha_v \frac{\ln(1 + \lambda \alpha_u \rho_v m_u)}{2} \;.
\end{equation*}
\begin{lemma}\label{lemma:properties_m_v^*}
For all $(m_u, \lambda) \in [0,\rho_u] \times (0,+\infty)$ there exists a unique $m_v^*(m_u,\lambda) \in [0,\rho_v]$ such that:
\begin{equation}
i(m_u, m_v^*(m_u,\lambda),\lambda) = \sup_{m_v \in [0, \rho_v]} i(m_u,m_v,\lambda)\;.
\end{equation}
Let $m_u(\lambda) \triangleq \rho_u \big(1-\frac{1}{1+\lambda \alpha_v \rho_u \rho_v}\big)$. The function $m_v^*(\cdot,\cdot)$ satisfies for $\forall (m_u, \lambda) \in \mathcal{D} \triangleq [0,\rho_u] \times (0,+\infty)$:
\begin{equation}
m_v^*(m_u, \lambda) =
\begin{cases}
\frac{m_u}{\lambda \alpha_v \rho_u (\rho_u-m_u)} \;\; \text{if} \quad 0 \leq m_u \leq m_u(\lambda)\;\;\,;\\
\qquad\; \rho_v \qquad\quad\: \text{if} \quad m_u(\lambda) < m_u \leq \rho_u \,.
\end{cases}
\end{equation}
It is continuous on $\mathcal{D}$ and continuously differentiable on $\mathcal{D}\setminus\{(m_u,\lambda) \in \mathcal{D}:m_u = m_u(\lambda)\}$. Finally,
\begin{equation}
\forall \lambda \in (0,+\infty),\forall m_u \in [0,m_u(\lambda)]:\frac{\partial i}{\partial m_v} \bigg\vert_{m_u,m_v^*(m_u,\lambda), \lambda} = 0 \;.
\end{equation}
\end{lemma}
\begin{IEEEproof}
Fix $(m_u, \lambda) \in \mathcal{D}$. Let $f: m_v \in [0,\rho_v] \mapsto i(m_u, m_v, \lambda)$. $f$ is continuously twice differentiable on $[0,\rho_v]$ and $f'(m_v) = \frac{\lambda \alpha_u \alpha_v}{2}\big(\frac{\rho_u}{1 + \lambda \alpha_v \rho_u m_v} - \rho_u + m_u \big)$, $f''(m_v) = -\frac{\lambda^2 \alpha_u \alpha_v^2 \rho_u^2}{2(1 + \lambda \alpha_v \rho_u m_v)^2}$.
We have:
\begin{equation*}
f'(m_v) = 0 \Leftrightarrow m_v = \frac{m_u}{\lambda \alpha_v \rho_u (\rho_u - m_u)}\;.
\end{equation*}
It is easy to check that this solution to $f'(m_v) = 0$ lies in $[0,\rho_v]$ if, and only if, $m_u \in [0,m_u(\lambda)]$ where $m_u(\lambda)$ is defined in the lemma. Besides, $f$ is strictly concave as $f''<0$.
Therefore, $f$ has a unique global maximizer that is given by the unique solution to $f'(m_v) = 0$ if $m_u \in [0,m_u(\lambda)]$ and is equal to $\rho_v$ if $m_u \in [m_u(\lambda), \rho_u]$.
The definition of $m_v^*(\cdot,\cdot)$ and its properties directly follows.
\end{IEEEproof}
\begin{lemma}\label{lemma:properties_h}
Let $h(\lambda) \triangleq \adjustlimits{\inf}_{m_u \in [0,\rho_u]} {\sup}_{m_v \in [0, \rho_v]} i(m_u,m_v,\lambda)$.
$h$ is continuously differentiable on $(0,+\infty)$ and satisfies for all $\lambda \in (0,+\infty)$:
\begin{align}
h(\lambda) &= i\big(m_u^*(\lambda) m_v^*(\lambda), \lambda \big) \;;\\
h'(\lambda) &= \frac{\alpha_u \alpha_v}{2}\big(\rho_u\rho_v - m_u^*(\lambda) m_v^*(\lambda)\big) \;;
\end{align}
where $(m_u^*(\lambda),m_v^*(\lambda))$ is the unique solution to the extremization that defines $h$:
\begin{align}
m_u^*(\lambda) &=
\begin{cases}
\qquad\quad\; 0\qquad\qquad\text{if}\; 0 < \lambda \leq \nicefrac{1}{\rho_u\rho_v \sqrt{\alpha_u\alpha_v}}\\
\frac{\lambda^2 \alpha_u \alpha_v \rho_v^2 \rho_u^2 - 1}{\lambda \alpha_u \rho_v (1 + \lambda \alpha_v \rho_v \rho_u)}\;\, \text{if}\; \lambda > \nicefrac{1}{\rho_u\rho_v \sqrt{\alpha_u\alpha_v}}
\end{cases}\;;\\
m_v^*(\lambda) &=
\begin{cases}
\qquad\quad\; 0\qquad\qquad\text{if}\; 0 < \lambda \leq \nicefrac{1}{\rho_u\rho_v \sqrt{\alpha_u\alpha_v}}\\
\frac{\lambda^2 \alpha_u \alpha_v \rho_v^2 \rho_u^2 - 1}{\lambda \alpha_v \rho_u (1 + \lambda \alpha_u \rho_v \rho_u)}\;\, \text{if}\; \lambda > \nicefrac{1}{\rho_u\rho_v \sqrt{\alpha_u\alpha_v}}
\end{cases}\;.
\end{align}
\end{lemma}
\begin{IEEEproof}
By Lemma~\ref{lemma:properties_m_v^*}, we have $h(\lambda) = \inf_{m_u \in [0,\rho_u]} g(m_u,\lambda)$ where $\forall \lambda \in (0,+\infty):g(m_u, \lambda) \triangleq i(m_u, m_v^*(m_u,\lambda),\lambda)$.
By continuity of $i(\cdot,\cdot,\lambda)$ and $m_v^*(\cdot,\lambda)$, $g(\cdot,\lambda)$ is continuous on $[0,\rho_u]$. Besides, $g(\cdot,\lambda)$ is increasing on $[m_u(\lambda),\rho_u]$ as $\forall m_u \in [m_u(\lambda),\rho_u]:$
$$
g(m_u,\lambda) = i(m_u, \rho_v,\lambda) = \alpha_u \frac{\ln(1+\lambda\alpha_v \rho_u \rho_v)}{2}
+ \alpha_v \frac{\ln(1 + \lambda \alpha_u \rho_v m_u)}{2} \;.
$$
Then, we can restrict the infimum to the interval $[0,m_u(\lambda)]$ in the definition of $h$: $h(\lambda) \triangleq \inf_{m_u \in [0,m_u(\lambda)]} g(m_u,\lambda)$.
For all $m_u \in [0,m_u(\lambda)]$:
\begin{align*}
g(m_u, \lambda) &= i\bigg(m_u, \frac{m_u}{\lambda \alpha_v \rho_u (\rho_u-m_u)},\lambda\bigg)\\
&= \frac{\alpha_u}{2 \rho_u}\big(\lambda \alpha_v \rho_v \rho_u^2 - m_u(1+\lambda \alpha_v \rho_v \rho_u)\big)
+ \frac{\alpha_u}{2}\ln\bigg(\frac{\rho_u}{\rho_u - m_u}\bigg)
+ \frac{\alpha_v}{2}\ln(1 + \lambda \alpha_u \rho_v m_u)\;;\\
\frac{\partial g}{\partial m_u}\bigg\vert_{m_u, \lambda}
&= -\frac{\alpha_u}{2 \rho_u}(1+\lambda \alpha_v \rho_v \rho_u)
+ \frac{\alpha_u}{2}\frac{1}{\rho_u - m_u}
+ \frac{\alpha_v}{2}\frac{\lambda \alpha_u \rho_v }{1 + \lambda \alpha_u \rho_v m_u}\\
&= \frac{\alpha_u \Big(-(1+\lambda \alpha_v \rho_v \rho_u)(\rho_u - m_u)(1 + \lambda \alpha_u \rho_v m_u)
	+ \rho_u(1 + \lambda \alpha_u \rho_v m_u)
	+ \lambda\alpha_v\rho_v \rho_u (\rho_u - m_u)\Big)}{2 \rho_u(\rho_u - m_u)(1 + \lambda \alpha_u \rho_v m_u)}\\
&= a(m_u, \lambda)\,
m_u\bigg(\frac{1-\lambda^2 \alpha_u \alpha_v \rho_v^2 \rho_u^2}{\lambda \alpha_u \rho_v (1 + \lambda \alpha_v \rho_v \rho_u)}
+ m_u\bigg)\;;
\end{align*}
where $a(m_u,\lambda) \triangleq \frac{\lambda \alpha_u^2 \rho_v (1 + \lambda \alpha_v \rho_v \rho_u)}{2 \rho_u(\rho_u - m_u)(1 + \lambda \alpha_u \rho_v m_u)}$.
Note that $\forall \lambda \in (0,+\infty),\forall m_u \in [0,m_u(\lambda)]:a(m_u,\lambda)>0$.
If $\lambda \leq \nicefrac{1}{\rho_u\rho_v \sqrt{\alpha_u\alpha_v}}$ then $0$ is the unique global minimizer of $g(\cdot,\lambda)$ on $[0, m_u(\lambda)]$.
Instead, if $\lambda > \nicefrac{1}{\rho_u\rho_v \sqrt{\alpha_u\alpha_v}}$, $\frac{\partial g}{\partial m_u}\big\vert_{m_u, \lambda}$ has a nonzero root given by:
\begin{equation*}
\frac{\lambda^2 \alpha_u \alpha_v \rho_v^2 \rho_u^2 - 1}{\lambda \alpha_u \rho_v (1 + \lambda \alpha_v \rho_v \rho_u)}
= m_u(\lambda)- \frac{1}{\lambda \alpha_u \rho_v (1 + \lambda \alpha_v \rho_v \rho_u)} \in (0,m_u(\lambda)) \;.
\end{equation*}
Then, it is easy to see that this root is the unique global minimizer of $g(\cdot,\lambda)$ on $[0, m_u(\lambda)]$.
We have just shown that
\begin{equation*}
m_u^*(\lambda) =
\begin{cases}
\qquad\quad\; 0\qquad\qquad\text{if}\; 0 < \lambda \leq \nicefrac{1}{\rho_u\rho_v \sqrt{\alpha_u\alpha_v}}\\
\frac{\lambda^2 \alpha_u \alpha_v \rho_v^2 \rho_u^2 - 1}{\lambda \alpha_u \rho_v (1 + \lambda \alpha_v \rho_v \rho_u)}\;\, \text{if}\; \lambda > \nicefrac{1}{\rho_u\rho_v \sqrt{\alpha_u\alpha_v}}
\end{cases}
\end{equation*}
is the unique global minimizer of $g(\cdot,\lambda)$ on $[0, m_u(\lambda)]$ (and, in fact, $[0,\rho_u]$) and $\forall \lambda \in (0,+\infty): \frac{\partial g}{\partial m_u}\big\vert_{m_u^*(\lambda), \lambda} = 0$.
Define
\begin{equation}
m_v^*(\lambda) \triangleq m_v^*(m_u^*(\lambda), \lambda)
\begin{cases}
\qquad\quad\; 0\qquad\qquad\text{if}\; 0 < \lambda \leq \nicefrac{1}{\rho_u\rho_v \sqrt{\alpha_u\alpha_v}}\\
\frac{\lambda^2 \alpha_u \alpha_v \rho_v^2 \rho_u^2 - 1}{\lambda \alpha_v \rho_u (1 + \lambda \alpha_u \rho_v \rho_u)}\;\, \text{if}\; \lambda > \nicefrac{1}{\rho_u\rho_v \sqrt{\alpha_u\alpha_v}}
\end{cases}\;.
\end{equation}
It follows from Lemma~\ref{lemma:properties_m_v^*} that $\forall \lambda \in (0,+\infty)$:
\begin{equation}\label{h_extremization_solved}
h(\lambda) = g(m_u^*(\lambda), \lambda) = i(m_u^*(\lambda), m_v^*(\lambda),\lambda)
\end{equation}
By Lemma~\ref{lemma:properties_m_v^*},
$\forall \lambda \in (0,+\infty):\frac{\partial i}{\partial m_v} \big\vert_{m_u^*(\lambda),m_v^*(\lambda), \lambda} = 0$
and
$m_v^*(\cdot,\lambda)$ is continuously differentiable on $[0,m_u(\lambda)]$ so (remember the definition of $g(m_u, \lambda)$ at the beginning of the proof):
\begin{equation*}
0 = \frac{\partial g}{\partial m_u}\bigg\vert_{m_u^*(\lambda), \lambda}
= \frac{\partial i}{\partial m_u}\bigg\vert_{m_u^*(\lambda), m_v^*(\lambda), \lambda}
+ \frac{\partial i}{\partial m_v}\bigg\vert_{m_u^*(\lambda), m_v^*(\lambda), \lambda} \cdot \frac{\partial m_v^*}{\partial m_u}\bigg\vert_{m_u^*(\lambda),\lambda}
= \frac{\partial i}{\partial m_u}\bigg\vert_{m_u^*(\lambda), m_v^*(\lambda), \lambda} \;.
\end{equation*}
All in all, we have shown that
\begin{equation}\label{di_dm_at_m*}
\forall \lambda \in (0,+\infty): \frac{\partial i}{\partial m_u} \bigg\vert_{m_u^*(\lambda),m_v^*(\lambda), \lambda} = \frac{\partial i}{\partial m_v} \bigg\vert_{m_u^*(\lambda),m_v^*(\lambda), \lambda} = 0 \;.
\end{equation}
Combining \eqref{h_extremization_solved} and the fact that $m_u^*, m_v^*$ are continuously differentiable on $(0,+\infty)\setminus\{\nicefrac{1}{\rho_u\rho_v \sqrt{\alpha_u\alpha_v}}\}$, we obtain that $h$ is continuously differentiable on $(0,+\infty)\setminus\{\nicefrac{1}{\rho_u\rho_v \sqrt{\alpha_u\alpha_v}}\}$ and for all $\lambda \in (0,+\infty)\setminus\{\nicefrac{1}{\rho_u\rho_v \sqrt{\alpha_u\alpha_v}}\}$:
\begin{align*}
h'(\lambda)
&=   \frac{\partial i}{\partial \lambda} \bigg\vert_{m_u^*(\lambda),m_v^*(\lambda), \lambda}
+ \frac{d m_u^*}{d \lambda}\bigg\vert_\lambda \cdot \frac{\partial i}{\partial m_u} \bigg\vert_{m_u^*(\lambda),m_v^*(\lambda), \lambda}
+ \frac{d m_v^*}{d \lambda}\bigg\vert_\lambda \cdot \frac{\partial i}{\partial m_v} \bigg\vert_{m_u^*(\lambda),m_v^*(\lambda), \lambda}\\
&= \frac{\partial i}{\partial \lambda} \bigg\vert_{m_u^*(\lambda),m_v^*(\lambda), \lambda}\\
&= \frac{\alpha_u \alpha_v}{2}\big(\rho_u-m_u^*(\lambda)\big)\big(\rho_v-m_v^*(\lambda)\big)
+ \frac{\alpha_u\alpha_v \rho_u m_v^*(\lambda)}{2(1+\lambda\alpha_v \rho_u m_v^*(\lambda))}
+ \frac{\alpha_u \alpha_v \rho_v m_u^*(\lambda)}{2(1 + \lambda \alpha_u \rho_v m_u^*(\lambda))}\\
&= \frac{\alpha_u \alpha_v}{2}\big(\rho_u-m_u^*(\lambda)\big)\big(\rho_v-m_v^*(\lambda)\big)
+ \frac{m_v^*(\lambda)}{\lambda}\bigg(\frac{\partial i}{\partial m_v} \bigg\vert_{m_u^*(\lambda),m_v^*(\lambda), \lambda} + \frac{\lambda\alpha_u\alpha_v}{2}\big(\rho_u-m_u^*(\lambda)\big)\bigg)\\
&\qquad\qquad\qquad\qquad\qquad\qquad\qquad\quad\;\,
+ \frac{m_u^*(\lambda)}{\lambda}\bigg(\frac{\partial i}{\partial m_u} \bigg\vert_{m_u^*(\lambda),m_v^*(\lambda), \lambda} + \frac{\lambda\alpha_u\alpha_v}{2}\big(\rho_v-m_v^*(\lambda)\big)\bigg)\\
&= \frac{\alpha_u \alpha_v}{2}\Big(\big(\rho_u-m_u^*(\lambda)\big)\big(\rho_v-m_v^*(\lambda)\big)
+ m_v^*(\lambda)\big(\rho_u-m_u^*(\lambda)\big) +m_u^*(\lambda)\big(\rho_v-m_v^*(\lambda)\big)\Big)\\
&= \frac{\alpha_u \alpha_v}{2}\big(\rho_u\rho_v-m_u^*(\lambda)m_v^*(\lambda)\big)\;.\\
\end{align*}
The second and second-to-last inequalities follow from \eqref{di_dm_at_m*}. Note that $\lambda \mapsto \frac{\alpha_u \alpha_v}{2}\big(\rho_u\rho_v-m_u^*(\lambda)m_v^*(\lambda)\big)$ is continuous at $\lambda = \nicefrac{1}{\rho_u\rho_v \sqrt{\alpha_u\alpha_v}}$. Therefore, $h$ is continuously differentiable on $(0,+\infty)$.
\end{IEEEproof}
We can now prove Theorem~\ref{theorem:MMSE} that we repeat here for reader's convenience:
\begin{theorem*}
Let $\lambda_{\mathrm{IT}} \triangleq \frac{1}{\rho_u\rho_v \sqrt{\alpha_u\alpha_v}}$.
For all $\lambda \in (0,+\infty)$, there is a unique solution to the extremization over $(m_u,m_v)$ on the right-hand side of \eqref{limit_mutual_information} given by:
\begin{equation*}
\big(m_u^*(\lambda),m_v^*(\lambda) \big)
=
\begin{cases}
	\qquad\qquad\qquad\;\;\: (0\,,0) &\text{if}\;\: 0 < \lambda \leq \lambda_{\mathrm{IT}}\\
	\Big(\frac{\lambda^2 \alpha_u \alpha_v \rho_v^2 \rho_u^2 - 1}{\lambda \alpha_u \rho_v (1 + \lambda \alpha_v \rho_v \rho_u)}
	,\frac{\lambda^2 \alpha_u \alpha_v \rho_v^2 \rho_u^2 - 1}{\lambda \alpha_v \rho_u (1 + \lambda \alpha_u \rho_v \rho_u)}\Big) &\text{if}\;\: \lambda > \lambda_{\mathrm{IT}}
\end{cases}\,.
\end{equation*}
	Then, $\mathrm{MMSE}_\lambda(\bU\bV^T \vert \bY)$ satisfies:
	\begin{equation*}
	\lim_{n \to +\infty} \mathrm{MMSE}_\lambda(\bU\bV^T \vert \bY) = \rho_u \rho_v - m_u^*(\lambda) m_v^*(\lambda)\;.
	\end{equation*}
	Hence, the asymptotic MMSE is less than $\rho_u \rho_v$ if, and only if, $\lambda > \lambda_{\mathrm{IT}}$.
\end{theorem*}
\begin{IEEEproof}
Let $n \in \mathbb{N}^*$. Define $h_n: \lambda \in (0,+\infty) \mapsto \frac{I(\bU,\bV;\bY)}{n}$ (the mutual information depends on $\lambda$ through the observation $\bY$).
The I-MMSE relation \cite{Guo2005_IMMSE} reads:
\begin{equation}\label{I-MMSE}
h_n'(\lambda)
= \frac{\partial}{\partial \lambda}\bigg(\frac{I(\bU, \bV ; \bY)}{n}\bigg)
= \frac{n_u}{n}\frac{n_v}{n} \frac{\mathrm{MMSE}_{\lambda}(\bU\bV^T \vert \bY)}{2}\;.
\end{equation}
$\lambda \mapsto \mathrm{MMSE}_{\lambda}(\bU\bV^T \vert \bY)$ is nondecreasing so $h_n$ is convex on $(0,+\infty)$.
By Lemma~\ref{lemma:properties_h}, $h: \lambda \mapsto \adjustlimits{\inf}_{m_u \in [0,\rho_u]} {\sup}_{m_v \in [0, \rho_v]} i(m_u,m_v,\lambda)$ is continuously differentiable on $(0,+\infty)$ and, by Theorem~\ref{theorem:limit_mutual_information}, it is the pointwise limit of the sequence of convex functions $(h_n)_{n \in \mathbb{N}^*}$. It follows:
\begin{equation*}
\forall \lambda \in (0,+\infty): \lim_{n \to +\infty} h_n^\prime(\lambda) = h'(\lambda) = \frac{\alpha_u \alpha_v}{2}\big(\rho_u\rho_v - m_u^*(\lambda) m_v^*(\lambda)\big) \;.
\end{equation*}
Combining the later with \eqref{I-MMSE} yields $\lim_{n \to +\infty} \mathrm{MMSE}_{\lambda}(\bU\bV^T \vert \bY) = \rho_u\rho_v - m_u^*(\lambda) m_v^*(\lambda)$.
The rest of the theorem follows from Lemma~\ref{lemma:properties_h}.
\end{IEEEproof}
\end{document}